\documentclass[%
 reprint,
superscriptaddress,
 amsmath,amssymb,onecolumn,
 aps
]{revtex4-2}

\pdfoutput=1
\usepackage[utf8]{inputenc} 
\usepackage[T1]{fontenc}    
\usepackage[dvipsnames]{xcolor}
\usepackage[colorlinks=true,citecolor=blue,
urlcolor=TealBlue]{hyperref}
\usepackage{url}            
\usepackage{booktabs}       
\usepackage{amsmath, amsthm, amssymb, amsfonts}
\usepackage{graphicx}
\usepackage{color}
\usepackage{array}

\usepackage{cleveref}
\crefname{figure}{Fig.}{Figs.}            
\crefname{equation}{Eq.}{Eqs.}            
\crefname{appendix}{Appendix}{Appendices} 
\crefname{section}{Section}{Sections}     
\crefname{enumi}{Step}{Steps}             
\crefname{table}{Table}{Tables}             

\creflabelformat{equation}{#2\textup{#1}#3} 
\usepackage{relsize} 
\usepackage{subcaption}

\usepackage[ruled,vlined]{algorithm2e}
\DeclareMathOperator\erf{erf}

\DeclareMathOperator\minmod{minmod}

\begin{document}

 \title{A hybrid level-set / embedded boundary method applied to
solidification-melt problems}
\author{A. Limare}
\affiliation{Sorbonne Universit\'e, CNRS, UMR 7190, Institut Jean Le Rond $\partial$'Alembert,  F-75005 \ Paris, France}
\affiliation{Laboratoire d'Hydrodynamique (LadHyX), UMR 7646 CNRS-Ecole Polytechnique, IP Paris, 91128 Palaiseau CEDEX, France}
\author{S. Popinet}
\affiliation{Sorbonne Universit\'e, CNRS, UMR 7190, Institut Jean Le Rond $\partial$'Alembert,  F-75005 \ Paris, France}
\author{C. Josserand}
\affiliation{Laboratoire d'Hydrodynamique (LadHyX), UMR 7646 CNRS-Ecole Polytechnique, IP Paris, 91128 Palaiseau CEDEX, France}
\author{Z. Xue}
\affiliation{Laboratoire d'Hydrodynamique (LadHyX), UMR 7646 CNRS-Ecole Polytechnique, IP Paris, 91128 Palaiseau CEDEX, France}
\author{A. Ghigo}
\affiliation{Department of Mathematics, British Columbia University, 1984 Mathematics Road, Vancouver, BC V6T 1Z4, Canada}
\date{}

\begin{abstract}

In this paper, we introduce a novel way to represent the interface for
two-phase flows with phase change. We combine a level-set method with a
Cartesian embedded boundary method and take advantage of both. This is part of
an effort to obtain a numerical strategy relying on Cartesian grids allowing
the simulation of complex boundaries with possible change of topology while
retaining a high-order representation of the gradients on the interface and
the capability of properly applying boundary conditions on the interface. This
leads to a two-fluid conservative second-order numerical method. The ability
of the method to correctly solve Stefan problems, onset dendrite growth with
and without anisotropy is demonstrated through a variety of test cases.
Finally, we take advantage of the two-fluid representation to model a
Rayleigh--B\'enard instability with a melting boundary.

\end{abstract}
\maketitle

\section{Introduction} 
Liquid--solid phase change (solidification or melting) is present in many
industrial processes, particularly in metallurgy~\cite{Chalmers1964} and 3-D
printing~\cite{Lewandowski2016}. Controlling ice formation and accretion is
also crucial in aeronautics with a recent increasing interest due to the
evolution of safety policies~\cite{Baumert2018,Villedieu2018}. More
generally, icing dynamics control a large number of important environmental
processes~\cite{Worster2000}, such as sea-ice
formation~\cite{Wettlaufer1997,Worster2006} or permafrost
thawing~\cite{Walvoord2016}. From an industrial point of view, reproducible
solidification processes which create complex geometries for solid materials with
isotropic properties at a low cost have been a goal pursued for
decades. Complex shape generation generally involves putting the
matter in a liquid state as an intermediary step before solidifying it, hence
the need to have a good knowledge of the process of solidification. This process
is difficult to study experimentally and often requires the use of intrusive or
sometimes destructive methods. Similarly, experimental studies on icing often
provide partial measurements only (surface temperature for instance) even
when they are made in controlled
conditions~\cite{Ghabache2016,schremb2016,Thievenaz2019,Thievenaz2020,ThievenazEPL,Monier2020}.

Numerical methods able to accurately simulate the process of solidification
and/or melting are thus of particular interest. Developing these methods is
especially challenging however, since melting and solidification processes
combine multiple difficulties. The first difficulty is classical and common to all free-boundary problems: how to accurately describe and follow the evolution of a complex boundary? This can be seen essentially as a geometric and kinematic problem and a broad range of methods have been proposed to solve it. The second difficulty concerns the dynamics of this motion (i.e. the relation between accelerations and forces) and requires the development of methods able to accurately couple the geometry of the boundaries with the underlying equations of motion. This coupling is clearly ``higher-order'' (in the sense of space/time derivatives) than the kinematic problem and thus more difficult to solve. A representative example is the approximation of surface tension terms which has been particularly challenging (see \cite{popinet2018numerical} for a review).

This coupling is especially difficult in the case of solidification/melting since the dynamics are driven almost entirely by singular terms on the boundary, such as temperature gradient jumps~\cite{Davis2006}. In the case of dendritic 
crystallisation the boundary topology can also become extremely complex and boundary-discontinuity difficulties can be compounded by the appearance of metastable states, for example in supercooled liquids.

A classical and accurate way to deal with partial differential equations with jumps is to use {\em boundary conforming} discretisation techniques combined with a Finite Volume Method (FVM) which ensures discrete and global mass conservation. In a boundary conforming framework, the mesh is constructed so that the edges of discretisation elements (for example triangles in 2D or tetrahedra in 3D) always coincide with the boundaries. Accurate jumps in the solutions can then be obtained by imposing the discrete boundary conditions directly on the edge of boundary elements. This allows in principle to design numerical schemes of arbitrary order of accuracy. The main limitation of these techniques is that they are inherently {\em
Lagrangian} i.e. they are most easily formulated in a Lagrangian frame of
reference and are thus in principle limited to small material deformations, such as occur for example in solid mechanics.  While techniques exist to overcome this limitation, such as Lagrangian-remapping \cite{loubere_subcell_2005}, they are usually complex and costly and still have difficulties dealing with complex topology changes such as merging and splitting.

This limitation of boundary-conforming techniques has led to the development of a broad range of methods able to couple general boundaries with the Eulerian framework more suitable to the discretisation of the equations of fluid motion. The issue then becomes: how to represent jumps/boundary conditions now that discrete boundaries do not coincide with real boundaries? The solution adopted by almost all methods to date is to approximate these (surface) jumps with localised {\em volumetric} terms which naturally fit within an Eulerian framework. This can be seen as replacing true Heaviside/Dirac functions with continuous/differentiable approximations and has a long history, dating back at least to the pioneering papers of Peskin \cite{peskin1972flow,peskin1977numerical}. 

A direct consequence of this approximation of discontinuous functions by differentiable approximations is that the resulting schemes can be at most first-order accurate spatially (by Godunov's theorem), in contrast with the boundary-conforming schemes mentioned earlier. This slow convergence is particularly problematic for applications which are mostly driven by interfacial terms, such as solidification and melting.

The goal of the present article is thus to lift this severe limitation and to present a Finite-Volume method able to deal with arbitrary boundary deformations, while conserving mass and preserving at least second-order spatial accuracy for the discretisation of boundary conditions and the overall solution.

\section{A brief review of existing schemes}

Non-boundary-conforming methods can can be classified in two families: front-tracking methods where one stores explicitly the position of the interface and front-capturing methods where the interface position is defined indirectly.

Juric and Tryggvason \cite{Juric1996b} for instance combined an explicit tracking of massless Lagrangian particles and an immersed boundary method. Another type of method based on the cellular automaton can also be used \cite{Gandin1994,Zhu2001,Zhu2002} often to study grain growth at the meso-scale. Reuther and Rettenmayr \cite{Reuther2014} simulated the dendritic solidification using an anisotropy-free meshless front-tracking method. However, the main drawbacks of these tracking methods are their difficulty to cope with change of topology and their complex extension to 3D.

In the second category, the interface is expressed implicitly using some auxiliary variables defined on every cell, for which values are ranging usually between zero and unity. Among others, one can cite the enthalpy method of Voller \cite{Voller2008} where the phase change occurs over a restricted temperature range and the solid-liquid interface is described as a mushy zone. Another family is the Volume of Fluid (VOF) method which ensures mass conservation \cite{hirt_arbitrary_1974,scardovelli_direct_1999}. Yet another widespread method is the phase field method \cite{Caginalp1986,karma1998quantitative,boettinger_phase-field_2002,plapp2010,Hester2020} which explicitly relies on a smooth, differentiable field representing phase transition. But, as discussed in the introduction, a large number of grid points in this transition zone is required for convergence.

The level-set method \cite{osher_fronts_1988,chen_simple_1997} is also a natural way to represent the interface which is simply a level set (usually the zero value) of a function defined in the calculation domain. Levet-set methods are well suited for modeling time-dependent, moving-boundary problems but also have their own specific drawbacks; they do not preserve mass/volume well in their original formulation, they introduce a smearing of the interface and reduce to low order accuracy regions where characteristics of the flow merge (i. e. caustic singularities). Furthermore, additional difficulties arise for the imposition of a flux jump condition on an interface and the associated construction of extension velocities. However, level-set methods are quite straighforward to implement, versatile enough to be combined with another method and their advantages and drawbacks, linked to the mathematical properties of the equations at play, have been studied quite thoroughly. 
Solutions have been found for applying an immersed boundary condition using a
finite-difference treatment for the variables, for instance the LS-STAG method
\cite{Cheny2010}, the Immersed Boundary Smooth Extension \cite{MacHuang2020}
or the Ghost Fluid method \cite{Fedkiw1999}. Note that all these methods can be shown to still rely on smooth approximations of Dirac/Heaviside functions \cite{popinet2018numerical}, and are thus only first-order accurate spatially.

On the other ahnd, cartesian embedded-boundary or \textit{cut-cell} methods have been extensively used for a large range of flows \cite{Popinet2003,hartmann2011strictly,berger2012progress}. They rely on a finite-volume discretization where cells are arbitrarily intersected by an embedded boundary. These methods show a second-order accuracy when applying immersed boundary condition and are conservative \cite{schwartz_cartesian_2006}. From an engineering point of view, this also greatly eases the mesh generation process. The main drawbacks of such methods are linked to grid irregularities in the cut regions which introduce local variations in truncation errors. This is all the more critical when the motion of the boundary is controlled by quantities calculated on the interface such as skin friction~\cite{Schneiders2013} or temperature gradients for phase change.

An important trend of the last two decades for numerical phase change models has been to create hybrid methods to compensate some of their shortcomings. For instance, phase change with VOF is especially hard since it has no built-in way of imposing Dirichlet conditions exactly on the interface, therefore it is often combined with other non-conservative methods which are able to impose a boundary condition on the interface. Sussman and Puckett \cite{Sussman2000} combined the VOF and level set methods; VOF ensures conservative properties whereas the level set method provides accurate geometric information such as normals and curvature. An extension of this method called CLSMOF was introduced in \cite{li2015incompressible}  with application to the freezing of supercooled droplets in~\cite{Vahab2016}. Recently, a hybrid VOF-IBM (Immersed Boundary Method) method has been developed for the simulation of freezing films and drops~\cite{legendre_2021}.

In the present article, this is precisely such a novel hybrid method that we introduce, by combining a level-set representation of the interface with a cut-cell method for the immersed boundary condition. By construction, this method is conservative and expected to a have a second-order accuracy.






\section{Principle of the method}

\begin{figure}[h!]
\centering
\def\svgwidth{0.5\textwidth}
    \centering
    \includegraphics[width = 0.5\textwidth]{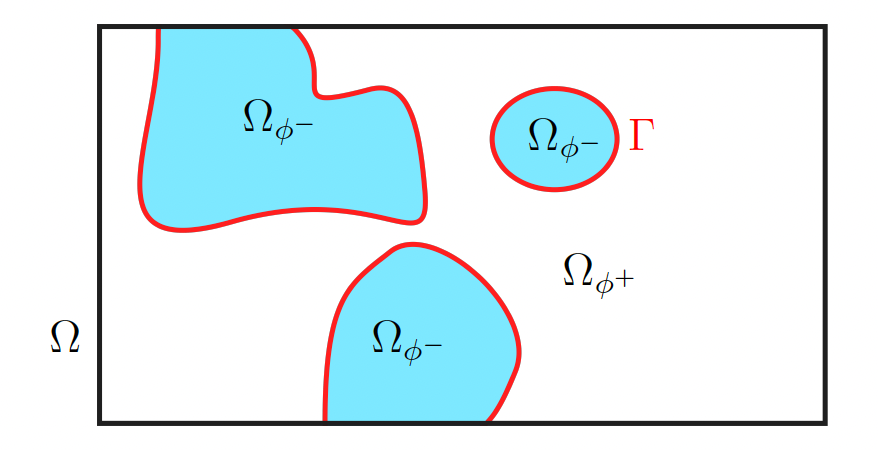}
\caption{Schematic view of domains used for calculation: in blue $\Omega_
{\phi^-}$, in white, $\Omega_{\phi^+}$, in red the interface $\Gamma$}
\label{fig:EmbedDomains}
\end{figure}

Most numerical methods take a ``one-fluid''
approach for multiphase flows, meaning that the computational domain
on which the numerical solver is applied contains both phases with a
more-or-less smooth change on the physical properties. Here, we
develop a two-fluid method, where each phase is described using its
own set of equations and variables. These two domains ($\Omega_
{\phi^+}$ and $\Omega_ {\phi^-}$ in Figure~\ref{fig:EmbedDomains}) are
coupled through the motion of the boundary $\Gamma$ and the associated
boundary/jump conditions. This approach has two main advantages which
are directly related to the similarities with boundary-conforming
Lagrangian methods: 1) The set of equations solved in one phase can be
different from those in the other phase (\textit{e.g.} a diffusion
equation in the solid and a Navier--Stokes equation coupled with
advection--diffusion in the fluid), 2) accurate boundary/jump
conditions can be imposed on the boundary. Specifically, the same
Dirichlet boundary condition is applied on the interface for both phases and the
(discontinuous) heat fluxes on the interface are calculated
independently for each phase with at least second-order accuracy using
finite-volume conservative numerical operators.

The boundary/interface $\Gamma$ is described using a levelset function $\phi$.
The domain outside of the interface is defined by $\Omega_ {\phi^+} =
\{\forall \boldsymbol{x} / \phi(\boldsymbol{x}) > 0\}$ and the inside of the
interface is defined in a similar manner $\Omega_ {\phi^-} = \{\forall
\boldsymbol{x} / \phi(\boldsymbol{x}) < 0\}$, both are subdomains of the
calculation domain $\Omega$. We depicted a possible situation on
\cref{fig:EmbedDomains}: in that case the blue domain represents $\Omega_{-}$
and is made of 3 disconnected subdomains. Let $T_ {S}$ and $T_{L}$ be
temperature fields defined respectively in $\Omega_{\phi^-}$ and
$\Omega_{\phi^+}$.

The temperature gradient jump then gives the velocity of the interface, which is the
starting point for the construction of a continuous extension velocity which
we will refer to as the phase change velocity $\boldsymbol{v}_{pc}$. The
assumption here is that the auxiliary field $\boldsymbol{v}_{pc}$ has a
meaning in both domains, not only on the interface $\Gamma$, thus allowing the
transport of the levelset function to solve the kinematic problem.

In the next section we present the physical model and the equations to be solved, which belong to the family of Stefan problems coupled with a velocity field~\cite{Gupta03}. The details of the coupling between the level-set method for the kinematic problem and the cut-cell technique \cite{johansen_cartesian_1998,schwartz_cartesian_2006} for the dynamic problem are given in section \ref{sec:numerical}. Finally in
\cref{sec:test_cases} we present several semi-analytical test cases and a more complex case by Favier \textit{et al.} \cite{Favier2019} where the equations solved
are simple diffusion in the solid and the Navier--Stokes equations in the
liquid.

\section{Physical model}\label{sec:equations}

The liquid--solid interface denoted $\Gamma$ separates two phases of a pure
material. Its position is determined by a prescribed temperature field at the
interface, that is not {\it a priori} constant and can depend on the
interface curvature and velocity, following for instance the Gibbs--Thomson relation. We consider that the solid domain cannot deform and that the liquid one obeys the
incompressible Navier--Stokes equations. Within this framework, the energy
equation simplifies into a diffusion equation in the solid domain and an
advection--diffusion equation in the liquid. The interface dynamics are
determined by the difference between the heat fluxes at the interface, following
the well-known Stefan equation. The solid and liquid parameters (diffusion
coefficients, viscosity and density in the fluid domain in particular) usually
depend on the temperature, but we will consider here constant values since
this paper focusses on the phase change dynamics, and the generalization to
smooth, temperature-dependent parameters does not bring additional numerical
challenges. We will only use the Boussinesq approximation for the Navier--Stokes
equation to model the Rayleigh--B\'enard thermal convective instability
during solidification. Our model includes the crucial physical
effects for solidification that are undercooling, crystalline anisotropy,
surface tension, and molecular kinetics. This allows us to treat problems with
supercooled fluids and study solidification fronts where instabilities occur
giving birth to dendrites and
fingering~\cite{ivantsov1947temperature,Mullins1964,Langer1980}. Our model
thus reduces to the following set of differential equations:
\begin{itemize}
  \item[-] for the temperature field $T_L({\bf x},t)$ in the fluid domain:
\begin{equation}
\rho_L C_{L}\left(\dfrac{\partial T_{L}}{\partial t}+{\bf u} \cdot \nabla T_L\right)= \nabla \cdot (\lambda_L \nabla T_L),
\label{eq:TL}
\end{equation}
\item[-] for the temperature field $T_S({\bf x},t)$ in the solid domain:
\begin{equation}
\rho_S C_{S}\dfrac{\partial T_{S}}{\partial t}= \nabla \cdot (\lambda_S \nabla T_S).
\label{eq:TS}
\end{equation}
\end{itemize}
where ${\bf u}({\bf x},t)$ is the velocity field (we consider that the
velocity vanishes in the solid domain). We denote with subscript $L$ and $S$
the coefficients related to the liquid and solid respectively and will use $i$
to denote either. $\rho_L$, $C_L$ and $\lambda_L$ ($\rho_S$, $C_S$ and
$\lambda_S$) are the liquid (solid) density, thermal capacity and thermal conductivity
respectively. The velocity field ${\bf u}({\bf x},t)$ in the fluid domain
obeys the incompressible Navier--Stokes equation that reads in its usual form:
\begin{eqnarray}
  \rho \left(\dfrac{\partial\bf u}{\partial t}+ {\bf u} \cdot \nabla ({\bf u})\right) &=& -\nabla p+\nabla \cdot ( 2 \mu {\bf D})+ \rho {\bf g}, \\
  \nabla\cdot{\bf u} &=& 0,
\label{eq:NS}
\end{eqnarray}
where ${\bf D}=\frac12 ({\bf \nabla}({\bf u})+^t{\bf \nabla}({\bf u}))$ is the
deformation tensor and ${\bf g}$ the acceleration of gravity. $\rho$ and $\mu$ are the density and
viscosity that can eventually depend on space through for instance the
temperature field. In this paper, we consider that the solver for the velocity
field in the liquid domain already exists and we will simply couple it with
the temperature field and solidification front dynamics. This coupling will be
used in the validation section \ref{sec:test_cases}, using the
Boussinesq approximation for the density variation with the temperature. This
set of differential equations needs to be complemented by the boundary
conditions at the interface $\Gamma$ describing the solidification front.
First, the temperature at the solidification front depends on the local
interface curvature and velocity, through the so-called Gibbs--Thomson
relation~\cite{Worster2000,Davis2006,Rappaz2011}
\begin{equation}
\forall \boldsymbol{x} \in \Gamma, \hspace{0.5 cm} T(\boldsymbol{x},t) = T_{m} -
\epsilon_
{\kappa}
\kappa - \epsilon_{v} v_{pc}
\label{eq:Tm}
\end{equation}
where $T_m$ is the melting temperature, $\kappa$ the local curvature of the
interface, $v_{pc}$  the local speed of the interface, $\epsilon_{v}$, the
molecular kinetic coefficient and $\epsilon_{\kappa}$ the surface tension
coefficient. Unless otherwise stated, these two coefficients will be taken as
constant in the present study.\\ Finally, the last equation couples the
thermal equations between the two domains (solid and liquid) stating that the
solidification front evolves through the balance between the heat flux at the
front, the so-called Stefan equation:
\begin{equation}
\rho_{S} L_H\boldsymbol{v}_{pc}\cdot \boldsymbol{n} =
\left[ \lambda\left.\nabla T\right|_{\Gamma} \right] \cdot \boldsymbol{n}
= \left(\lambda_{L}\left.\nabla T_{L}\right|_{\Gamma} -
\lambda_{S}\left.\nabla T_{S}\right|_{\Gamma}\right)\cdot \boldsymbol{n}
\label{eq:Stefan}
\end{equation}
where $L_H$ is the latent heat and $\boldsymbol{n}$ the normal to the
interface from solid to liquid. The velocity of the interface
$\boldsymbol{v_{pc}}$ is thus related through the latent heat to the jump in
the heat flux (and therefore in general to the temperature gradient) across
the interface.

Finally, a dimensionless version of this set of equations will be used, introducing reduced temperature, geometrical length and time scales. They usually lead to a dimensionless Stefan number, that compares thermal diffusion and latent heat. Its definition might depend on the specificity of the problem (geometry, boundary conditions). For the sake of simplicity, we will consider later on that the density of the liquid and the solid are the same ($\rho_L=\rho_S=\rho$): although it is not true in general (for the ice/water phase change, we have $\rho_S/\rho_L \sim 0.9$ for instance), it is not a crucial ingredient for the numerics~\cite{legendre_2021}.
Considering a domain of size $L_0$ and defining thus a time scale $\tau=L_0^2/D_S$) using the thermal diffusion coefficient in the solid ($D_S=\lambda_S/(\rho_S C_S)$), we obtain the following set of dimensionless equation (defining also $D_L=\lambda_S/(\rho_L C_L)$):

\begin{equation}
\frac{\partial \theta_L}{\partial t}+{\bf u} \cdot \nabla \theta_L=\frac{D_L}{D_S} \Delta \theta_L,
\label{ad:TL}
\end{equation}
\begin{equation}
\frac{\partial \theta_S}{\partial t}=\Delta \theta_S,  \,\,\, {\rm and}
\label{ad:TS}
\end{equation}
\begin{equation}
\boldsymbol{v}_{pc}\cdot \boldsymbol{n} = {\rm St} \left( \frac{\lambda_{L}}{\lambda_S} \left.\nabla \theta_{L}\right|_{\Gamma} -
\left.\nabla \theta_{S}\right|_{\Gamma}\right)\cdot \boldsymbol{n}.
\label{ad:Stefan}
\end{equation}
We have introduced a reduced temperature, defined using a temperature $T_1$, coming in general from the boundary conditions and thus depending on the specific problem to investigate, leading typically:
$$ \theta_{L,S}= \frac{T-T_m}{T_1-T_m}=\frac{T-T_m}{\Delta T}.$$
Here $\Delta T=T_1-T_m$ is supposed to be positive, leading to the following definition of the Stefan number:

$$ {\rm St}=\frac{C_S \Delta T}{L_H}.$$
The Stefan number thus quantifies the ratio between the available heat in the system with the latent heat. In the following we will in general use this set of dimensionless equations, noting the dimensionless temperature $T$ instead of $\theta$ by simplicity.

\section{Numerical method}\label{sec:numerical}
The goal of this paper is to present a numerical method able to solve
accurately the thermal equations (\ref{eq:TL}, \ref{eq:TS}, \ref{eq:Tm} and
\ref{eq:Stefan}), that will be coupled with an existing solver for the fluid
equation (\ref{eq:NS}). The method will be implemented in the free software
\textit{Basilisk}\cite{basilisk}. We use a novel approach for the numerical modelling of the
interface by combining a level-set function with an embedded boundary
(cut-cell) treatment for the fluxes. This means that for interfacial cells we
store two different values for the temperature fields  in order to correctly
compute the temperature gradient in each phase. We will first describe the
choices made for the level-set function and the Cartesian embedded-boundary
and then explain how we combined both of these approaches to obtain a
consistent numerical description of the physical situation.

\subsection{Global algorithm}\label{subsec:generalAlgo}
Our method can be summarized as:
\begin{enumerate}
    \item Calculate the phase change velocity on the interface
    $\left.\boldsymbol{v_{pc}}\right|_{\Gamma}$ \label{item:initial_value}
    \item Extend or reconstruct a continous phase change velocity field ${\boldsymbol v_{pc}}$ in the
    vicinity of the interface from the $\left.\boldsymbol{v_{pc}}\right|_{\Gamma}$ value
    \label{item:extension}
    \item Advect the level-set function and recalculate the volume and face
    fractions \label{item:LSAdvect}
    \item Redistance the level-set function
    \item Initialize fields of newly emerged cells\label{item:initialization}
    \item Apply the appropriate solver for each independent phase
    \item Perform mesh adaptation\label{item:mesh_adapt}
\end{enumerate}

Key points that will be further detailed are
\cref{item:initial_value,item:extension} which combine the level-set
representation of the interface for the reconstruction of a continuous field
with the calculation of the gradients on the interface relying on the embedded
boundary representation of the interface, \cref{item:initialization} that is
critical for the global accuracy of the method and \cref{item:mesh_adapt}
which allows efficient calculations.

\subsection{The level-set method}
\label{subsec:LS}
The level-set is a method initially designed to study the motion by a velocity
field $\boldsymbol{v}$ of an interface $\Gamma$ of codimension 1 that bounds
several open regions $\Omega$ (possibly connected)\cite{gibou2018review}. The main idea is to use a
function $\phi$ sufficiently smooth (Lipschitz continuous for instance) and
define the interface as the 0-level-set of $\phi$:
\begin{equation}
\forall \boldsymbol{x}\in\Gamma\hspace{0.1cm} , 
\hspace{0.2cm} \phi(\boldsymbol{x},t) = 0
\end{equation}
and the equation of motion of the level-set function is:
\begin{equation}
\dfrac{\partial \phi}{\partial t} + \boldsymbol{v}\nabla \phi = 0
\label{eq:LS_adv}
\end{equation}
where $\boldsymbol{v}$ is the desired velocity on the interface. The level-set
method has multiple advantages, the main one for our calculations being built-in
topological regularization that deals easily with merging and pinching off,
and allows robust calculation of geometric properties.\\

Even though the 0-level-set will be advected with the correct velocity, $\phi$
will no longer be a distance function and can become irregular after several
timesteps. Because the values of the level-set function in the vicinity of the
0-level-set are used to reconstruct a velocity field (see \cref{subsec:speed_recons}), it
hinders this reconstruction process, hence the need to correct the values of
the level-set to get $|\nabla \phi | = 1$. One way is to iterate on the
following Hamilton-Jacobi equation \cite{Sussman1994}:
\begin{equation}
\left\{\begin{aligned}
\phi_\tau + sign(\phi^{0}) \left(\left| \nabla \phi\right| - 1 \right) &= 0\\ 
\phi(x,0) &= \phi^0(x)
\end{aligned}
\right.
\label{eq:LS_redist}
\end{equation}
where $\tau$ is a fictitious time and $\phi^0$ is the value of $\phi$ at the
beginning of the redistancing process. Numerous methods for reinitialization
exist, see \cite{Solomenko2017} for a comparative study or the recent work of
Chiodi and Desjardins \cite{Chiodi2017}. We took the method of Min \& Gibou
\cite{Min2007} with corrections by Min \cite{Min2010}, derived from the method
of Russo \& Smereka \cite{russo_remark_2000}. In order to preserve the mass
and the position of 0-level-set $\phi_0$, the idea is to include the initial
interface location in the stencils of the discretized spatial derivatives. We
have recalled the procedure which can be naturally extended to 3D in
\ref{sec:Redistancing}.\\


\subsection{Embedded Boundary (cut-cells)}
\label{subsec:embed}
\begin{figure}[h!]
    \centering
    \includegraphics[width = 0.5\textwidth]{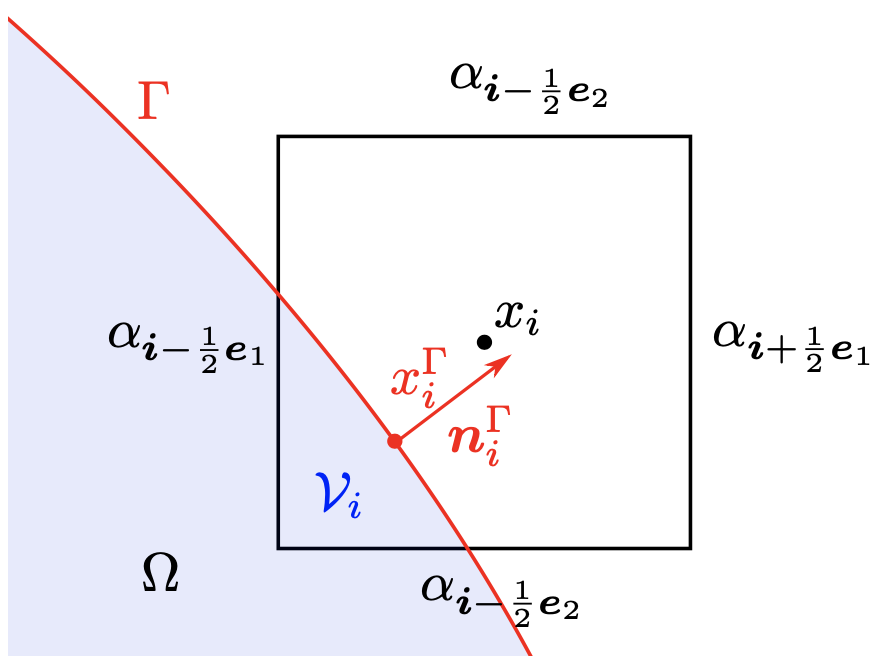}
    \caption{Cut-cell with its initial center of gravity $x_{i}$ outside the
  calculation
  domain $\Omega$
  }
  \label{fig:cut_cell_johansen}
\end{figure}

The 0-level-set of the distance function defines the
interface $\Gamma$ between the two phases. This level-set function is used as
an input to modify control volumes in a finite-volume manner. In this section,
we draw the main lines of the embedded boundary method as defined in
\cite{johansen_cartesian_1998,schwartz_cartesian_2006} and introduce the
notations used hereafter. The main idea is to consider a domain $\Omega$ with
a general boundary $\Gamma $ and embed this domain in a regular
Cartesian grid with $\Delta x$ the grid spacing and $\boldsymbol{e}_{d}$
the unit vector in the $d$ direction. The intersection of each cell with
$\Omega$ gives a collection of irregular cells as shown on \cref{fig:cut_cell_johansen}. The vertex-centered levelset field can then be used to obtain the volume fractions $\mathcal{V}_{i}$ defined as
\begin{equation}
  |V_{i}| = \mathcal{V}_{i}\Delta x^{D}
  \label{eq:volumefrac}
\end{equation}
with $|V_{i}|$ the volume of a cell and $D$ the dimension of the problem (2
or 3). Similarly, the face fractions
$\alpha_{\boldsymbol{i}\pm\frac{1}{2}\boldsymbol{e_{d}}}$ are defined as
\begin{equation}
|A_{\boldsymbol{i}\pm\frac{1}{2}\boldsymbol{e}_{d}}| =  \alpha_{\boldsymbol{i}\pm\frac{1}{2}\boldsymbol{e_{d}}} \Delta x^{D-1}
\label{eq:facefrac}
\end{equation}
where $|A_{\boldsymbol{i}\pm\frac{1}{2}\boldsymbol{e}_{d}}|$ is the surface of the face. The details of the calculation of the volume fractions and face fractions can be found on Basilisk's website (\href{http://basilisk.fr/src/embed.h}{http://basilisk.fr/src/embed.h}) and is a direct adaptation of the Johansen \& Collela's \cite{johansen_cartesian_1998}. These fractions give access to collection of piecewise linear segments, whose centroids $x_{i}^{\Gamma}$ and
normals $\boldsymbol{n}_{i}^{\Gamma}$ are defined as
\begin{equation}
\boldsymbol{x}_{i}^{\Gamma} =
\dfrac{1}{|\mathcal{A}_{i}^{\Gamma}|}\int_{\mathcal{A}_ {i}^{\Gamma}}
\boldsymbol{x}d\mathcal{A}.
\end{equation}
\begin{equation}
\boldsymbol{n}_{i}^{\Gamma} = 
\dfrac{1}{|\mathcal{A}_{i}^{\Gamma}|}\int_{\mathcal{A}_ {i}^{\Gamma}}
\boldsymbol{n}_{i}^{\Gamma}d\mathcal{A}.
\end{equation}
which, in turn form the basis for the construction of conservative, high-order discretization operators,  especially the divergence operator $\nabla$:
\begin{equation}
\left.\begin{aligned}
\nabla \cdot \overrightarrow{F} &\approx \dfrac{1}{|V_{i}|}\int_{V_{i}} \nabla
\cdot \overrightarrow{F} dV = \dfrac{1}{|V_{i}|}\int_{\partial V_{i}}
\overrightarrow{F}\cdot \boldsymbol{n}dA\\
& \approx \left[\left(\sum_{\pm=+,-}\sum_{d=1}^{D}\pm \alpha_ {\boldsymbol{i}\pm
\frac{1}{2}\boldsymbol{e_{d}}}F^{d}(\boldsymbol{x}_{i\pm \frac{1}{2}
\boldsymbol{e_{d}}})\right) + \alpha_{i}^{\Gamma}\boldsymbol{n}_{i}^
{\Gamma}\cdot
\overrightarrow{F}(\boldsymbol{x}_{i}^{\Gamma})\right]
\end{aligned}\right..
\end{equation}

\begin{figure}[h!]
    \centering
    \includegraphics[width = 0.5\textwidth]{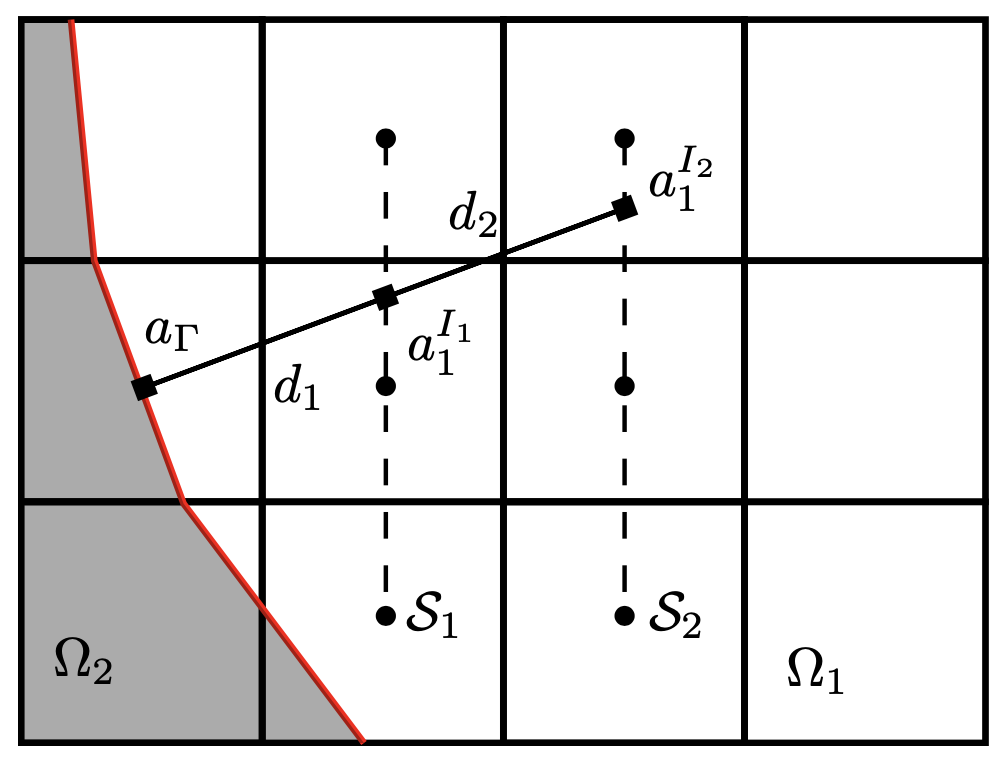}
\caption{Gradients calculation}
\label{fig:Embed1}
\end{figure}

The key point here is that the equations solved on the separate subdomains $\Omega_{\phi^\pm}$ bounded by $\Gamma$ are sufficiently smooth (since they do not include the interfacial discontinuities) and can therefore be
extended to the domain made of the cells that have a non-zero volume fraction,
even when the cells have their original center outside of the calculation domain. We depict on \cref{fig:cut_cell_johansen} a possible configuration where a cell has its center $x_{i}$ outside of the domain $\Omega$ bounded by $\Gamma$.
This leads to a conservative, finite-volume methodology which is at least second-order
accurate. The main advantage of Cartesian grids embedded boundary over
structured or unstructured grid methods is simpler grid generation. The
underlying regular grid also allows the use of simpler data structures and
numerical methods over a majority of the domain. Accuracy is maintained at the
boundaries using an algorithm detailed in
\cite{johansen_cartesian_1998}. For each partially covered cell or
interfacial cell (cells for which $0<\mathcal{V}_{i}<1$), the flux through the
boundary, which is the crucial ingredient of \cref{eq:Stefan}, is calculated using only values from other cells.\\
The gradient of a variable $a_{1}$ defined only in one phase of the
calculation domain $\Omega_{1}$, phase 1 on \cref{fig:Embed1}, on the
embedded boundary $\left. \nabla a_{1}\right|_{\Gamma}$ in the direction of $\boldsymbol{n_{\Gamma}}$ is calculated as
\begin{equation}
\left. \nabla a_{1}\right|_{\Gamma} = 
\dfrac{1}{d_{2}-  d_{1}}\left[ \dfrac{d_{2}}{d_{1}}(a_{\Gamma}- a^{I_{1}}_{1}) 
- \dfrac{d_{1}}{d_{2}}(a_{\Gamma}- a^{I_{2}}_{1})\right]
\label{eq:grad_embed}
\end{equation}
where $a_{1}^{I_{1,2}}$ are quadratically interpolated values of $a_{1}$ on each segment $\mathcal{S}_{1,2}$ and 
$a_{\Gamma}$ is the imposed Dirichlet boundary condition on the interface
(here the Gibbs--Thomson relation) which is the same for both phases. The stencil used
for the interpolation follows the procedure described in~\cite{schwartz_cartesian_2006},
 depending on the normal of the interface $\boldsymbol{n}=
\{n_{l},l=1,\dots,d\}$, here the dimension $d=2$, the segments
$\mathcal{S}_{1}$ and $\mathcal{S}_{2}$ are chosen to be normal to
$\boldsymbol{e_{k}}$ where $\{k : n_{k}\geq n_{l}, l=1,2\}$.\\
The calculation of the gradient of $\left.\nabla a_{2}\right|_{\Gamma}$ in
$\Omega_{2}$, is done by using a second variable $a_{2}$ defined only in the
second domain, volume fractions and face fractions of the second calculation
domain can be deduced from the ones previously calculated, they are just the
complementary to 1.

The two scalar variables controlling the motion of the interface are the
temperature variables denoted $T_{L}, T_{S}$, the temperature
fields for the liquid and the solid respectively. For each temperature field, we apply
the embedded boundary method and build the associated discretization operators
independently using the appropriate volume and face metrics. The interface
temperature $T_{\Gamma}$  given by \cref{eq:Tm} is required to calculate the
temperature gradients. We thus need two quantities:
\begin{itemize}
  \item[-] the curvature $\kappa$, calculated using the height
  function\footnote{With a level-set function, it would be tempting to use the classical
  relation for the calculation of the curvature
  \begin{equation}
  \kappa = \dfrac{\phi_y^2\phi_{xx} - 2\phi_x\phi_y\phi_{xy}+\phi_x^2\phi_{yy}}
  {\left| \phi_x^2 + \phi_y^2\right|^{3/2}},
  \end{equation}
  preliminary tests showed no difference between the 2 methods. For a more detailed comparison of the influence of curvature on dendritic growth, one may refer to \cite{lopez2013two}.} as in
  Popinet~\cite{popinet_accurate_2009},
  \item[-] the phase change velocity $v_{pc}$, we assume that we have
  previously calculated the velocity of the interface and that this velocity
  remains constant during a timestep,
\end{itemize}
the details of these calculations will be discussed in \cref{subsec:time_disc_specificity}. This
yields the temperature gradients $\left.\nabla T_{L}\right|_{\Gamma}$ and
$\left. \nabla T_{S}\right|_{\Gamma}$ with second-order accuracy. Recalling
that the velocity of the interface is defined as the jump in the normal
direction of the interface of the gradient of the temperature fields
$L_{H}\boldsymbol{v_{pc}|_{\Gamma}} =
\left[\lambda\left.\nabla T\right|_{\Gamma} \right] \cdot \boldsymbol{n}$, we also
expect second-order accuracy on the velocity of the interface.

\subsection{Speed reconstruction off the interface}
\label{subsec:speed_recons}
\begin{figure}[h!]
    \centering
    \includegraphics[width = 0.5\textwidth]{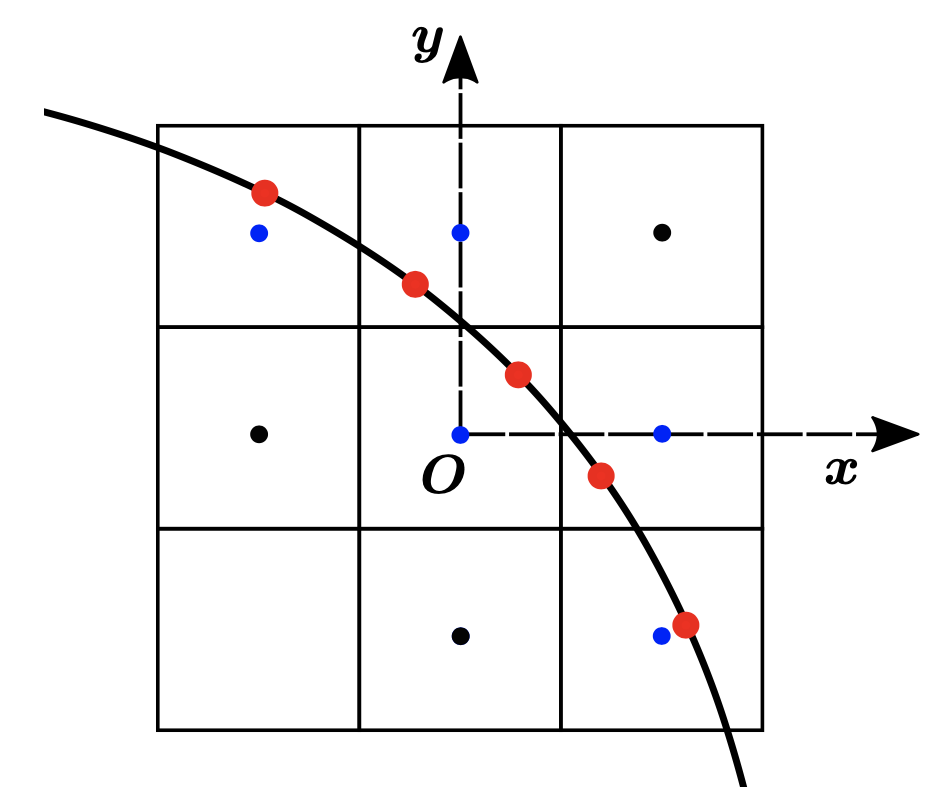}
\caption{Cell categories for the reconstruction of $ v_{pc} $ off the
interface. In blue: interfacial cells, in red: interface centroids, in black:
non-interfacial cells}
\label{fig:Recons1}
\end{figure}

To use the phase change velocity as the velocity in the level-set advection
equation, \cref{eq:LS_adv}, we have to
build a continuous velocity in the vicinity of the interface. In this section we
now describe how we rely on our level-set function for this process, starting
from the discrete velocity defined
only on the interface $\left.v_{pc}\right|_{\Gamma}$ that we have previously calculated.
We follow the approach of Peng \textit{et al.} \cite{peng_pde-based_1999}, and
solve an additional PDE so that
$\boldsymbol{v_{pc}}$ is constant along a curve normal to $\Gamma (\tau)$
\begin{equation}
\dfrac{\partial \boldsymbol{v_{pc}}}{\partial \tau}  + \delta S(\phi) 
\boldsymbol{n_
{\phi}}. \nabla \boldsymbol{v_{pc}} = 0\label{eq:speedrecons}
\end{equation} 
where $\delta$ is equal to 0 in interfacial cells and 1 elsewhere, $S(\phi)$
is the sign function and the vector $\boldsymbol{n_ {\phi}} =
\frac{\nabla\phi}{|\nabla \phi|}$ is normal to the isovalues of the level-set
function. The velocity reconstruction process is initialized by setting $
\boldsymbol{v_{pc}}(x,\tau=0) =\left.
\boldsymbol{v_{pc}}\right|_{\Gamma}$ in the interfacial cells
(blue cells on \cref{fig:Recons1} with the value of interfacial centroids in
red) and $0$ elsewhere. We want to highlight here that this field
reconstruction step couples the embedded boundary representation of the
interface using the interface centroids and the level-set method which gives
the normal to the interface $\boldsymbol{n_{\phi}}$.\\

\noindent The reconstruction is divided in two steps on
which we iterate until we reach convergence:
\begin{enumerate}
    \item A few iterations of \cref{eq:speedrecons} are performed, typically
    $2D$ with $D$ the dimension of the problem,
    such that the velocity converges in the vicinity of the interface. Note that the definition of $\delta$ ensures that only
    non-interfacial cell values are updated, the initial value in interfacial
    cells acts as a source term.
    \item The value of the velocity in interfacial cells is modified as
\begin{equation}
\widetilde{\boldsymbol{v_{pc}}} = \boldsymbol{v_{pc}} + \dfrac{\epsilon}{
\alpha_{ij}}
\label{speed_correction}
\end{equation}
where $\epsilon$ is the error between the interpolation of $\boldsymbol{v_
{pc}}$ on the boundary face centroids:
\begin{equation}
\forall \boldsymbol{x}_{i,\Gamma} \in \text{facets},
\hspace{0.5cm}\mathcal{L}(\boldsymbol{v_{pc}} (\boldsymbol{x}_{i,\Gamma}))  = 
\left. \boldsymbol{v_{pc}}\right|_{\Gamma} + \epsilon 
\label{eq:interp_recons}
\end{equation}
with $\mathcal{L}$ a biquadratic interpolation operator using Lagrange polynomials on a standard $3\times3$ stencil such that 
\begin{equation}
  \mathcal{L}(q (\boldsymbol{x})) = \sum_{i,j=-1\dots1}\alpha_{ij} q_{ij}
  \label{eq:biquadInterp}
\end{equation}
is the interpolation of a cell-centered field at $\boldsymbol{x}$.
\end{enumerate}
The aim of the second step is to correct the approximation error done at the initialization of this reconstruction process, by setting the value in the cell centers equal to the value in the centroids.

\begin{figure}[h!]
    \centering
    \includegraphics[width = 0.5\textwidth]{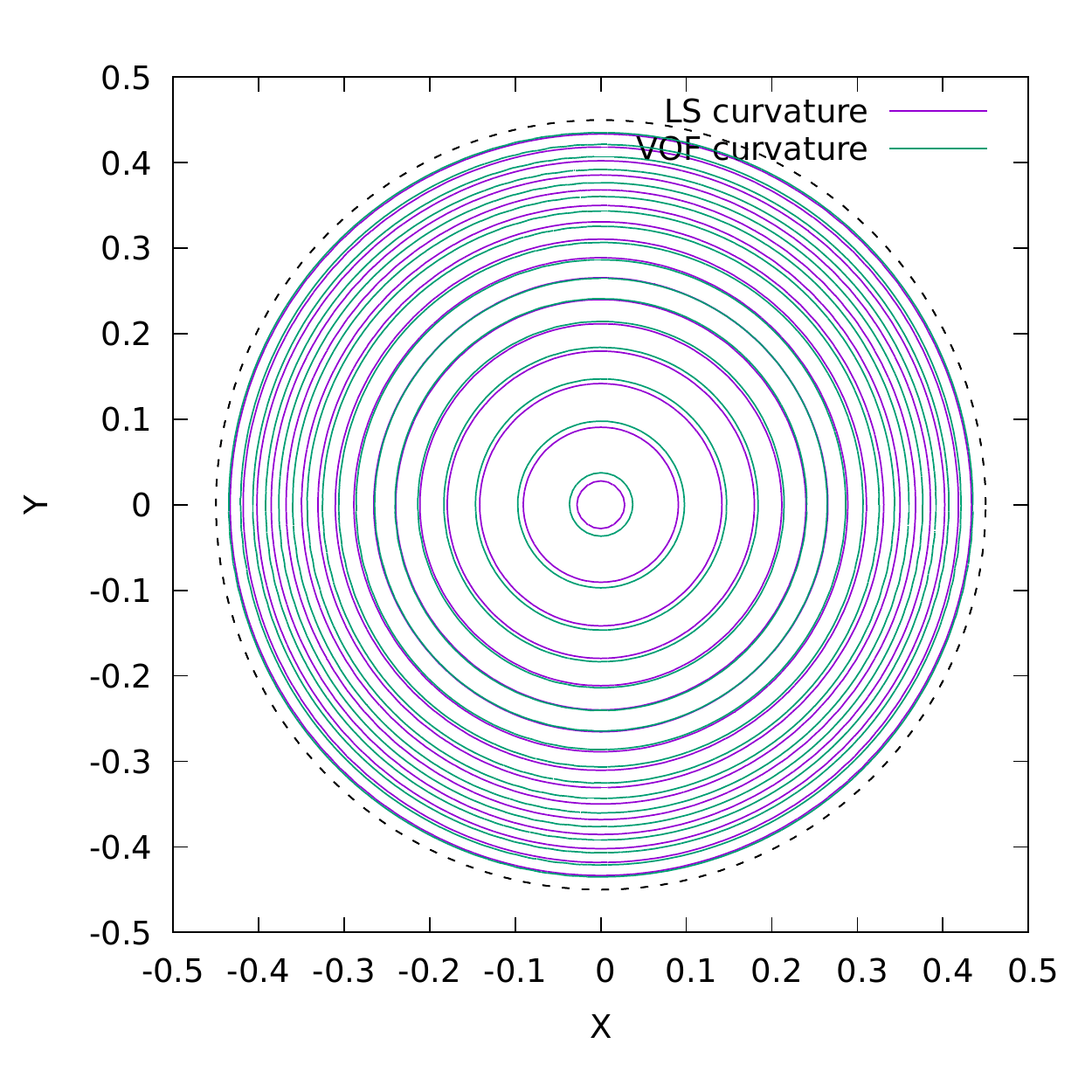}
    \caption{Shrinking circle, dashed line: initial interface}\label{fig:shrinking_circle}
\end{figure}

\begin{figure}[h!]
\centering
\begin{minipage}[h!]{.25\linewidth}
\centering
\includegraphics[width=0.99\textwidth]{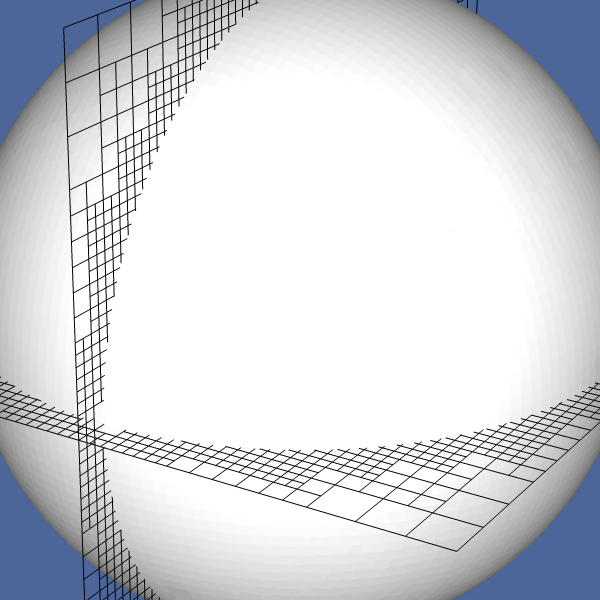}
\subcaption{$t = 0$}
\end{minipage}
\begin{minipage}[h!]{.25\linewidth}
\centering
\includegraphics[width=0.99\textwidth]{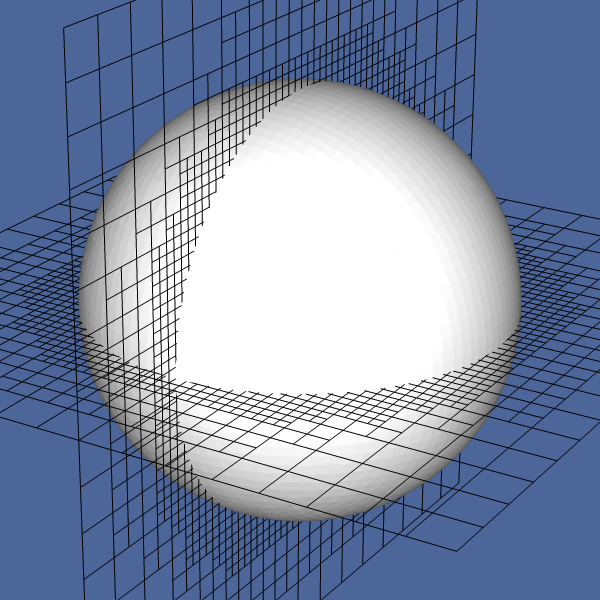}
\subcaption{$t = 120 \Delta t$}
\end{minipage}
\begin{minipage}[h!]{.25\linewidth}
\centering
\includegraphics[width=0.99\textwidth]{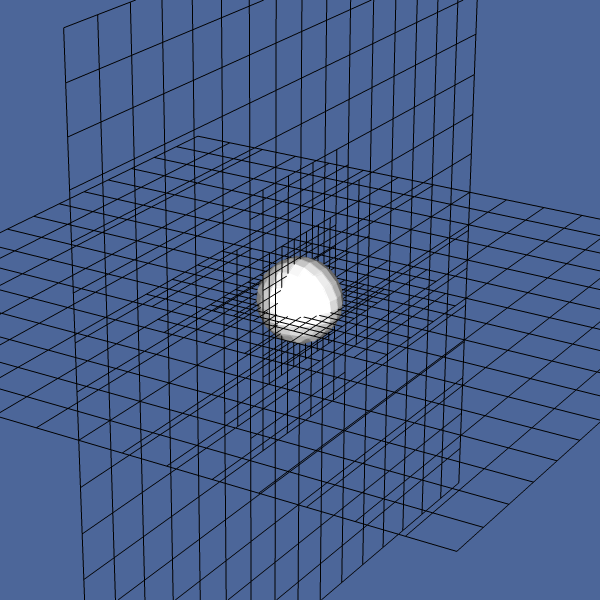}
\subcaption{$t = 240 \Delta t$}
\end{minipage}
\caption{Shrinking sphere, $\boldsymbol{v} = \kappa\boldsymbol{n}$ imposed on the
interface centroids}
\label{fig:growing_sphere}
\end{figure}

We show on \cref{fig:shrinking_circle} an application of this reconstruction
method, the initial calculation domain is $[-0.5:0.5]\times[-0.5:0.5]$. The initial
interface is a circle of diameter $9/10^{\text{th}}$. The initial grid
maximal resolution is $128\times128$. We set $\left.\boldsymbol{v_
{pc}}\right|_{\Gamma} = -\boldsymbol{n}\kappa $ in interfacial cells, with $\kappa$ the local curvature. The same test case has been run  using the curvature calculated with the height function and the level-set function. The
extension method for the velocity is applied with a CFL number of $0.3$. The level-set
function is then advected with the continuous velocity. At the end of each
iteration the level-set function is reinitialized. The interface is
output every 60 iterations and remains circular for both methods, as shown on \cref{fig:shrinking_circle}.
This demonstrates both the robustness and the accuracy of the method without any
additional regularization. We obtain similar results in 3D, as
shown in
\cref{fig:growing_sphere}.

\subsection{Embedded boundary motion: timestep constraint, emerging cells
scalar field initialization and truncation error variations}
\label{subsec:time_disc_specificity}

\begin{figure}[h!]
    \centering
    \includegraphics[width = 0.5\textwidth]{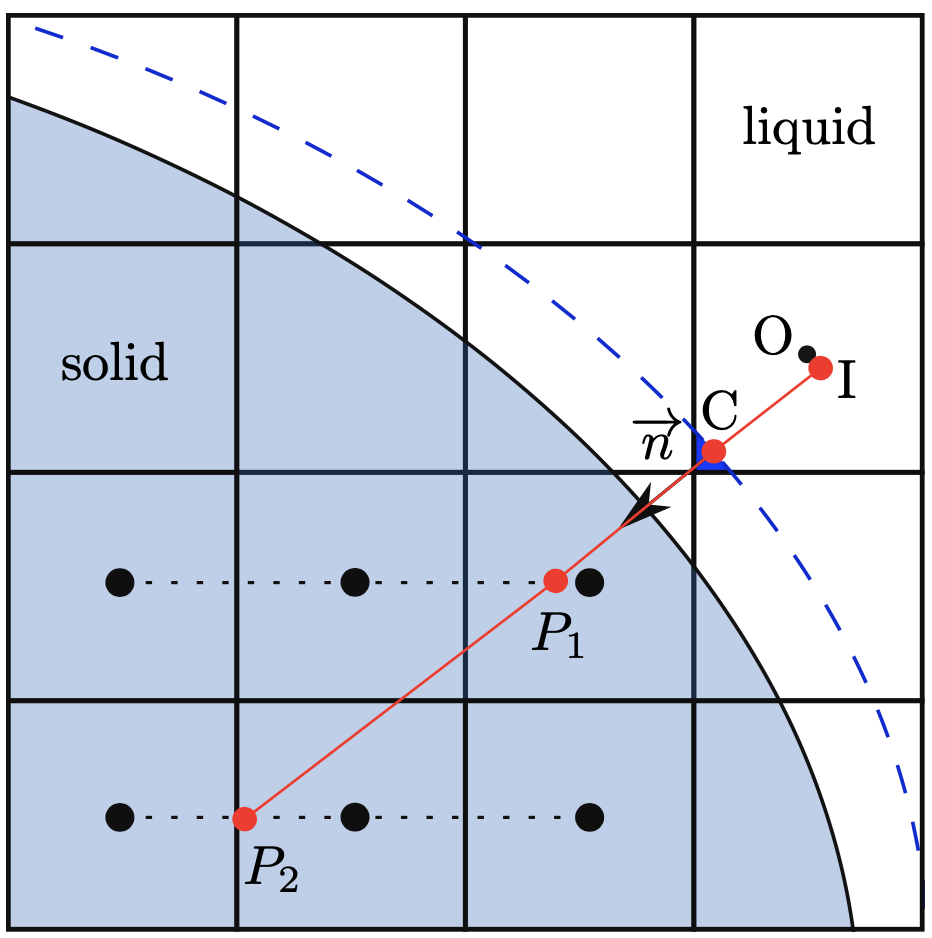}
\caption{Case of a new interfacial cell, in dark blue. Solid line:
$t^{n-1}$ and dashed line: $t^{n}$}
\label{fig:new_interf}
\end{figure}

In the case of a moving interface $\Gamma$, the two considered domains are
functions of time $\Omega_{+} = \Omega_{+}(t)$ and $\Omega_{-} =
\Omega_{-}(t)$. Cells that were non-interfacial can become interfacial. In these
cells an initialization technique for the undefined fields is required. We
show a typical case of a moving boundary on \cref{fig:new_interf} where the
interface $\Gamma$ at instant $t^{n-1}$ and instant $t^{n}$ are displayed with a
solid line and a dashed line respectively. The blue cell for which the solid
temperature was undefined at instant $t^{n-1}$ becomes an interfacial cell after
displacement and the solid temperature field needs to be initialized in this
cell. 

More generally, we tag new interfacial cells if they verify 2 conditions:
\begin{equation*}
\left\{\begin{aligned}
\mathcal{V}_{ij}^{n-1} (1-\mathcal{V}_{ij}^{n-1}) &=&0\\
\mathcal{V}^{n}_{ij} (1-\mathcal{V}^{n}_{ij}) &\neq& 0
\end{aligned}\right..
\end{equation*}
where we denote $\mathcal{V}^{n-1}$ the volume fraction at
the previous timestep and $\mathcal{V}^{n}$ at the current one, this assertion
is tested at \cref{item:LSAdvect} of the global algorithm. In order to
initialize the fields of emerging cells, we took an approach
similar to what is already done for gradient calculations. We detail here the initialization procedure only for the temperature. The
Dirichlet boundary condition on the interface $T^{n}_{\Gamma}$ is obtained from the geometric properties
of the interface after advection and an interpolation on the interface of
the previously reconstructed phase change velocity field
$\mathcal{L}(\boldsymbol{v_{pc}}^{(n)} (\boldsymbol{x}_{i,\Gamma}))$:
\begin{equation}
T_{L,\Gamma}^{n} = T_{S,\Gamma}^{n} = T_{\Gamma}^{n} = T_{m}
-\epsilon_{\kappa}\kappa^{(n)}
- \epsilon_{v}\mathcal{L}(\boldsymbol{v_{pc}}^{(n-1)} (\boldsymbol{x}_
{i,\Gamma})).
\label{eq:new_gibbs}
\end{equation}
The temperature field is interpolated, if possible, at two points along the
direction of the normal to the embedded boundary at $P_{1}$ and $P_{2}$. The
temperature at the centroid $C$ is given by the Gibbs--Thomson relation. A
quadratic interpolation of the values $P_ {1}, P_{2}, C$ gives the value at
$I$, the orthogonal projection of $O$ onto the line
$(C,\overrightarrow{n})$ where $\overrightarrow{n}$ is the normal to the
interface at $C$. We neglect tangential variations of the
temperature. Similar procedure can be devised for other variables (velocity, pressure) according to their boundary conditions on the interface.\\

Without any timestep constraint, some cells that were completely uncovered
might become completely covered. This means that a cell could undergo a
complete phase change during one timestep. Therefore, the
following constraint is applied:
\begin{equation}
\Delta t < \dfrac{\Delta x}{|\boldsymbol{v_{pc}}|}.
\end{equation}
The volume fractions $\mathcal{V}_{i}$ and face fractions $\alpha_{
\boldsymbol{i}\pm\frac12\boldsymbol{e}_{d}}$ are considered constant during
one timestep. We solve a fixed-boundary problem at each timestep, see Eqs.
(18-19) in Schwartz \textit{et al.}~\cite{schwartz_cartesian_2006}. Thus, our
numerical scheme for the displacement is only first-order accurate in time. This is a
strong approximation because in diffusion-driven cases, the motion timescale of
the interface is comparable to the diffusion timescale. Future work should
focus on using a better approximation of the position of the interface during
a timestep.

\begin{figure}[h!]
    \centering
    \includegraphics[width = 0.4\textwidth]{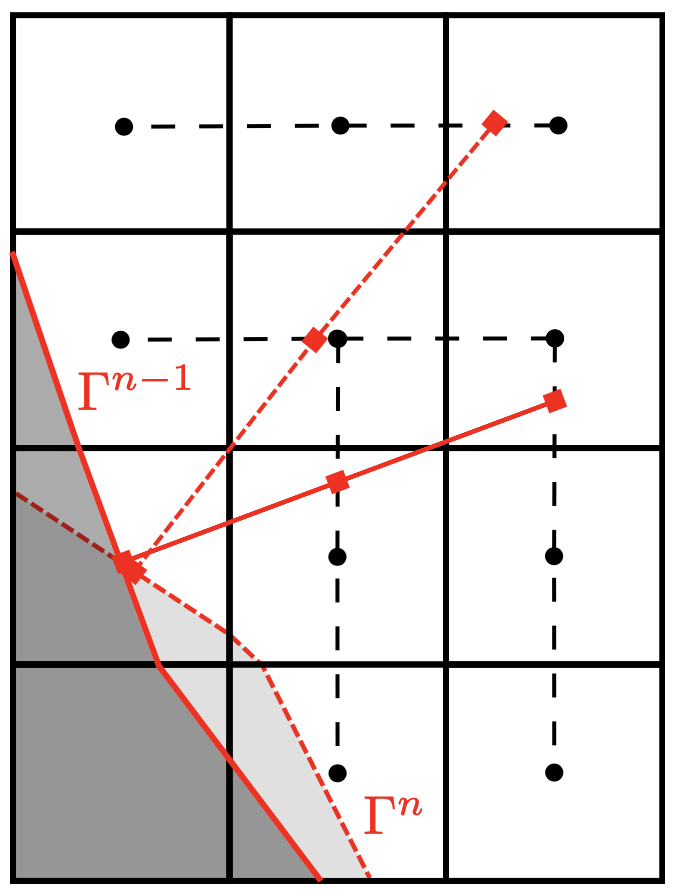}
    \caption{Stencil jump with interface motion, solid red line: $\Gamma^
    {n-1}$, dashed: $\Gamma^{n}$}\label{fig:truncation_temporal}
\end{figure}

Another issue related to Cartesian cut-cell methods with moving
embedded boundaries is the oscillation of fluxes calculated on the boundary
due to varying truncation errors. In cut cells, the discretization stencils
are offset and can vary abruptly as the interface moves, simply because the interface becomes cut/uncut. Another situation is shown on \cref{fig:truncation_temporal}, the red solid line shows the interface $\Gamma^{n-1}$ at
instant $t^{n-1}$, and the solid line is $\Gamma^{n}$. As one can see the
stencil used for the calculation of the fluxes varies and can introduce spurious
oscillations. Therefore, obtaining smoothly varying values of fluxes on the
interface is an active field of research 
\cite{Schneiders2013,berger_ode-based_2017}. Because
these variations of the truncation error interact with the motion of the
embedded boundary they can quickly deteriorate the quality of the solution.

\subsection{Mesh adaptation}
\begin{figure}[h!]
    \centering
    \includegraphics[width = 0.5\textwidth]{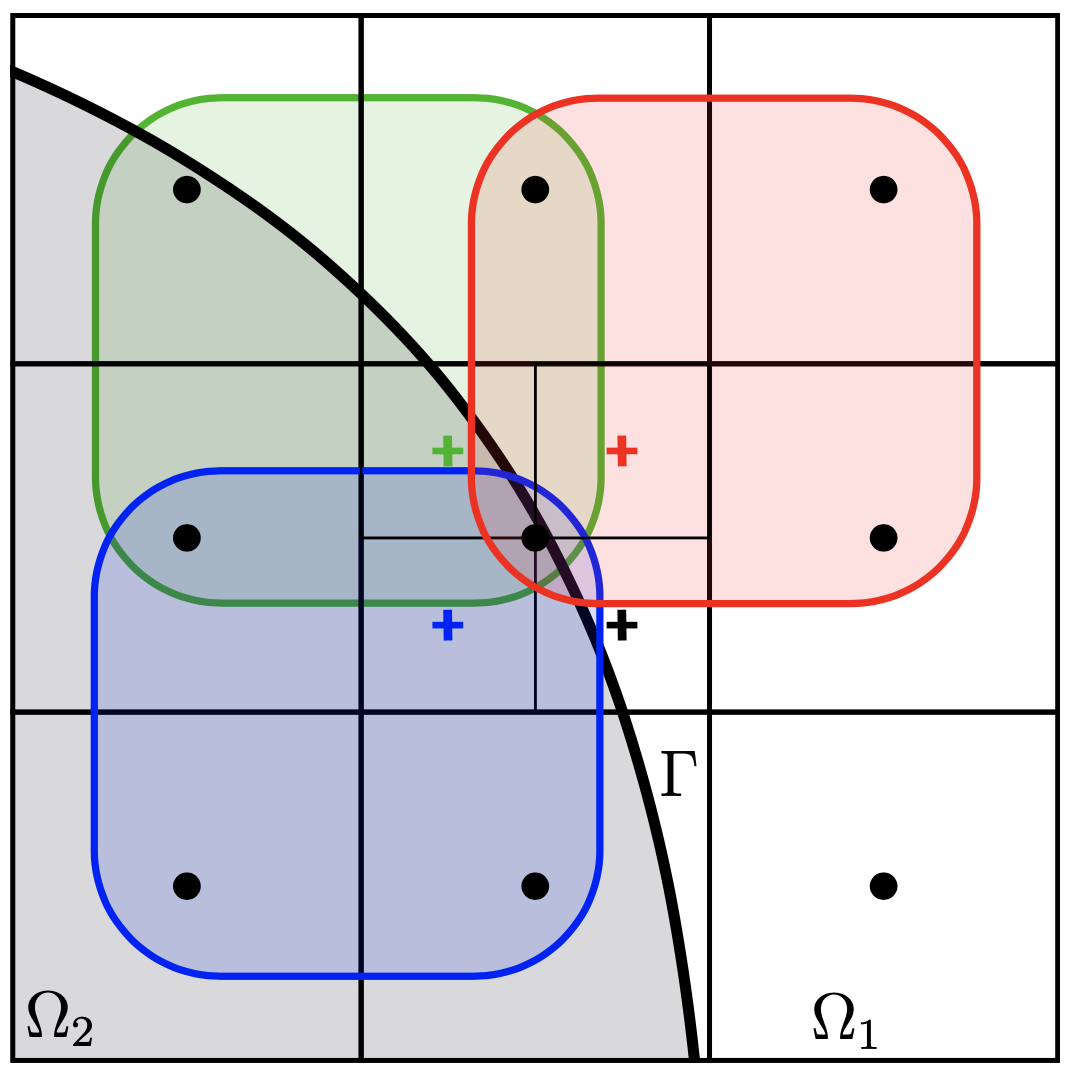}
  \caption{Mesh adaptation with an embedded boundary. Possible configurations for the prolongation operator.}
  \label{fig:metric_mesh_adaptation}
\end{figure}

The \textit{Basilisk} library has mesh adaptation capabilities. It uses
quad/octrees with a 2-1 balancing rule, see
\cite{popinet2015quadtree,van2018towards}. When one cell is refined, projected
values onto new cells are typically calculated via a bi/trilinear interpolation of the
field on the coarser level. In the presence of an embedded boundary, specific
refinement and coarsening functions have been written to take the face and volume fractions into account. We present on \cref{fig:metric_mesh_adaptation} a typical $3\times3$
stencil where the central cell is cut by the interface $\Gamma$. We detail
here the case of 3 of the children cells for a scalar field $a_{1}$ defined
only in the domain $\Omega_{1}$, we will refer to the sub-cells by the color
of their cell center:
\begin{itemize}
   \item red cell: a standard bilinear interpolation can be applied. The associated stencil
for interpolation, denoted by a red rectangle, only contains cells that are either partially covered or fully uncovered,
  \item green cell: only three cells on the coarser level are accessible,
  therefore a triangular interpolation is used,
  \item blue cell: this cell is completely covered, therefore it does not need to be initialized.
 \end{itemize} 
The same characterization is done in phase $\Omega_{2}$
simultaneously for the prolongation of $a_{2}$, which gives
\begin{itemize}
   \item red cell: this cell is partially covered, on the coarser grid, the
   diagonal cell is completely covered, therefore it is initialized with the
   value of its parent,
   \item green cell: triangular interpolation,
   \item blue cell: bilinear interpolation.
 \end{itemize}
\noindent The phase change velocity is used as an adaptation criterion
in our simulations in combination with the other ``standard'' adaptation
criteria, namely the temperature, the velocity in the liquid phase and the
volume fraction.

\section{Test \& validation cases}\label{sec:test_cases}

In this section we present eight different numerical configurations 
to demonstrate and characterize the ability of our method to obtain accurate
solutions, all of code associated for running those simulations is available in \href{http://basilisk.fr/sandbox/alimare/README}{A. Limare's sandbox}. The first case involves a planar interface and validates the
accuracy of the method without the initialization procedure for emerging
cells. The second is similar but this times validates the scalar field
initialization procedure. The third one tests the stability of the method in
2D and in particular the speed reconstruction method with a standard case known as the Franks's spheres. The fourth one
illustrates the ability of our method to capture instabilities and in particular
the formation of dendrites; it also checks and quantifies the accuracy of our
method. The fifth shows the compatibility of our method with anisotropy in the
Gibbs--Thomson condition (\href{http://basilisk.fr/sandbox/alimare/cube_sixfold.c}{cube\_sixfold.c}). The sixth is exactly the same simulation 
 one but in 3D and shows the formation of dendrites in 3D (\href{http://basilisk.fr/sandbox/alimare/crystal_growth3D.c}{crystal\_growth3D.c}). The seventh one makes a comparison of the tip velocity calculated using our method with a linear solvability theory. The
last simulation is taken from \cite{Favier2019} and combines two different
solvers: in the liquid and solid phase we solve the coupled diffusion
equations for the temperature as for the previous cases, but in the liquid, we
now solve in addition the Navier-Stokes equations allowing for fluid motion (\href{http://basilisk.fr/sandbox/alimare/Favier_Ra-Be.c}{Favier\_Ra-Be.c}).
In particular, we recover some of the main results of their study, which are
the existence of a critical Rayleigh number for the instability and the
variation with time of the associated wavelength because of the melting of the
solid boundary. Except when it is explicitly stated, we will consider for the
example the dimensionless set of equations
(\cref{ad:TL,ad:TS,ad:Stefan}) considering the fluid at rest (${\bf
u}=0$ everywhere) and taking the thermal ratios unity,
$D_L/D_S=\lambda_L/\lambda_S=1$ (recall that we have already taken
$\rho_L=\rho_S=\rho$).

\subsection{Solidifying domain}
This test case is borrowed from \cite{chen_simple_1997}. It is a simple Stefan
problem where an initially planar interface translates at constant velocity.
The interface is located at $t=0$ at the position $x=0$ with the initial temperature field:
\begin{equation}
T_0(x) = \left\{\begin{aligned}
-1 + e^{-V\cdot x}& , x > 0\\
0               & , x \leq 0\\
\end{aligned}\right.
\end{equation}
It is easy to show that the solution of the diffusion equation for the temperature field and the Stefan condition for the interface leads to the translation of the
planar interface at constant velocity $V=1$. Indeed, the temperature field evolves as
\begin{equation}
T(x,t) = \left\{\begin{aligned}
-1 + e^{-V(x-Vt)}& , x > Vt\\
0               & , x \leq Vt\\
\end{aligned}\right.
\end{equation}
which gives the steadily moving planar interface, whose equation is
\begin{equation}
\Gamma (t) = \{x = Vt, y=s\}\text{ , } s\in \mathfrak{R}.
\end{equation}
We perform an error analysis of our method by studying the error on the initial
phase change velocity, see \cref{tab:planar1D_initial_gradient}. It shows a second-order accuracy
on the initial temperature gradient jump and a classical error analysis on the final temperature field shows also a
second-order accuracy, see 
\cref{tab:planar1D}, with a slight drop in the order of accuracy for low
resolution. Note that we do not really validate the accuracy of the
initialization procedure of the temperature in the solid here, because the newly
solid cells only need to be initialized with $T=0$. Finally we did an error
analysis with a fixed timestep $\Delta t = 1\times10^{-6}$ and a fixed number of
400 iterations: the results are presented in \cref{tab:planar1D_fixed} and also show 
second-order asymptotic accuracy.\\

\begin{table}[h!]
    \centering
    \begin{tabular}{|c|c|c|}
        \hline
        Grid    & $L_{1}$-error & order\\\hline
        32$^2$  & 3.18e-04  & --    \\
        64$^2$  &  8.04e-05 & 1.98 \\
        128$^2$ &  2.02e-05 & 1.99 \\
        256$^2$ &  5.07e-06 &    2 \\\hline
    \end{tabular}
    \caption{Convergence of the initial temperature gradient when refining grid
    size}
    \label{tab:planar1D_initial_gradient}
\end{table}

\begin{table}[h!]
    \centering
    \begin{tabular}{|c|c|c|c|c|c|}
        \hline
        Grid    & Timestep            & $L_{1}$-error & order   & $L_{\infty}$-error & order\\\hline
        32$^2$  & $1.6\times 10^{-3}$ &  1.59e-4     & --      &  5.31e-4          & --   \\
        64$^2$  & $4  \times 10^{-4}$ &  6.52e-05     & 1.28    &  0.000252          & 1.07 \\
        128$^2$ & $1  \times 10^{-4}$ &  1.55e-05     & 2.07    &  6.46e-05          & 1.96 \\
        256$^2$ & $2.5\times 10^{-5}$ &  4.06e-06     & 1.93    &  1.63e-05          & 1.99 \\ \hline
    \end{tabular}
    \caption{Convergence of the temperature field when refining grid size and time step}
    \label{tab:planar1D}
\end{table}

\begin{table}[h!]
	\centering
	\begin{tabular}{|c|c|c|c|c|}
		\hline
		Grid    & $L_{1}$-error & order & $L_{\infty}$-error & order\\\hline
		32$^2$  &  1.51e-05     & --    &  1.67e-4           & -- \\
		64$^2$  &  5.52e-06     & 1.45  &  8.78e-05          & 0.92   \\
		128$^2$ &  1.4e-06      & 1.97  &  2.28e-05          & 1.94  \\
		256$^2$ &  3.32e-07     & 2.07  &  4.73e-06          & 2.26  \\ \hline
	\end{tabular}
	\caption{Convergence of the temperature field with grid refinement, fixed timestep = $10^{-6}$}
	\label{tab:planar1D_fixed}
\end{table}


\subsection{Planar interface with an expanding liquid domain}

\begin{figure}[h!]
    \centering
    \includegraphics[width = 0.4\textwidth]{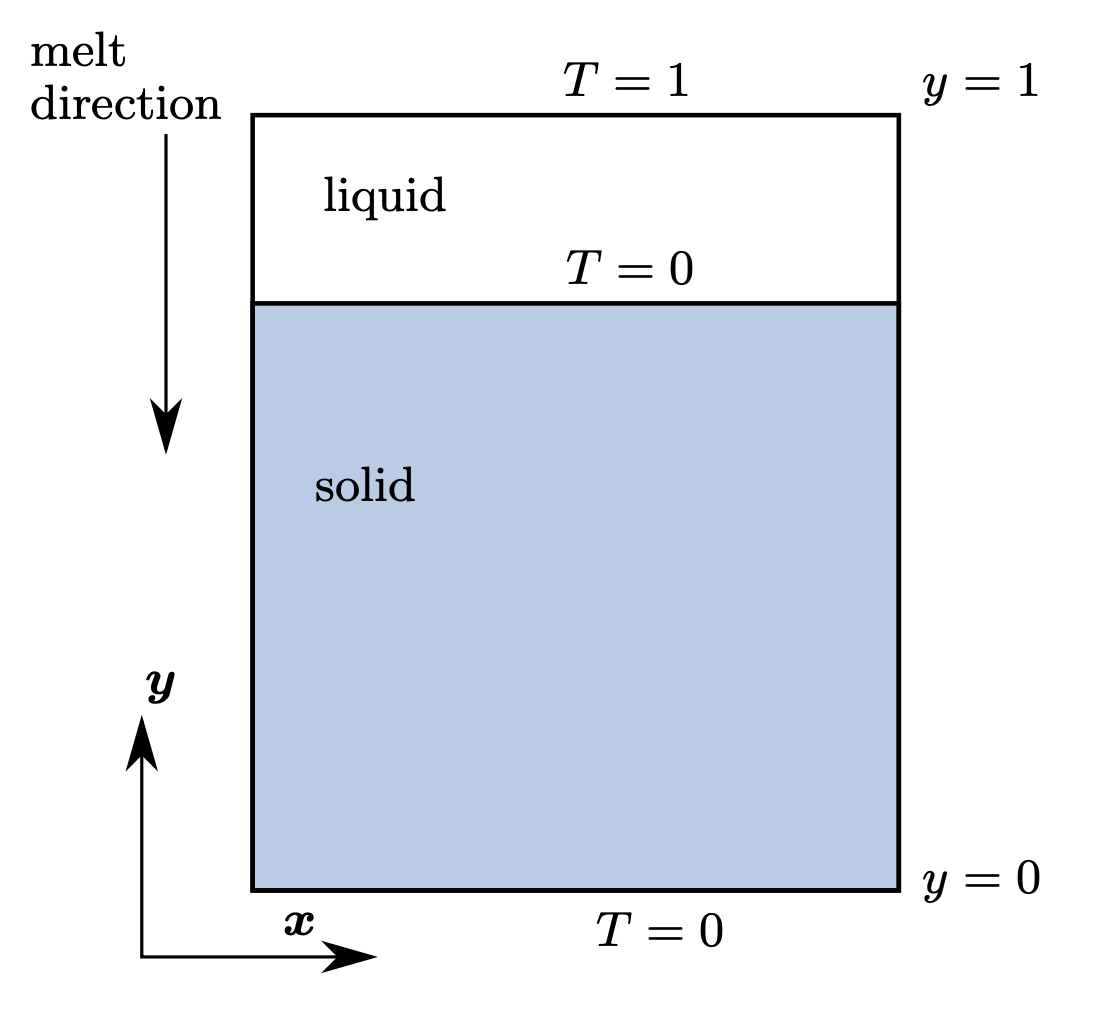}
  \caption{Scheme of the configuration studied}\label{fig:planar_expand_scheme}
\end{figure}

This case tests the diffusion of two tracers separated by an embedded boundary (taken from Crank
\cite{Crank1987}).  It corresponds to the melting of an ice layer by imposing a warm temperature condition $T_1$ at the top boundary and the melting one $T_m$ at the bottom. The Stefan number $St = \frac{C(T_1-T_{m} }{L_{H}}$ for our simulation is $2.85$ and the dimensionless temperature ($(T-T_m)/(T_1-T_m)$ still denoted 
$T$ in the dimensionless equation) is $1$ at the top ($y=1$) and $0$ at the bottom ($y=0$) as
shown in \cref{fig:planar_expand_scheme}. The initial temperature in the liquid is
\begin{equation*}
T_{L}(x,y,t_0) = 1- \frac{\erf(\dfrac{1-y}{2\sqrt{t_0}})}{\erf(\lambda)}
\end{equation*}
where $\lambda = 0.9$. The temperature in the solid is $T_{S} = 0$. With these
initial conditions, the interface position as a function of time is given by
\begin{equation}
y(t) = 1- 2\lambda \sqrt(t).
\end{equation}
We start the simulation at $t=t_{0} = 0.03$ such that there are at least two
full cells above the interface in order to have a correct approximation of the
gradients for \cref{eq:grad_embed}, in the liquid phase. Notice that the
initialization method of the temperature field in newly liquid cells is thus tested for this set up. Error
plots in \cref{tab:planar1D_fixed} shows convergence of the $L_{1}-norm$
with an order of accuracy slightly above $2$. The results on the $L_
{\infty}$-error also show the expected order of accuracy for low resolution and a drop at $256\times256$ which requires further investigation. 

\begin{table}[h!]
  \centering
  \begin{tabular}{|c|c|c|c|c|c|}
    \hline
    Grid    & Timestep & $L_{1}$-error & order & $L_{\infty}$-error & order\\\hline
    32$^2$  & 1.e-2   &  1.97e-03  & --    &  4.86e-03   & --     \\
    64$^2$  & 2.5e-3  &  3.80e-04  & 2.37  &  6.40e-04   & 2.93\\
    128$^2$ & 6.25e-4 &  8.31e-05  & 2.19  &  1.41e-04   & 2.17\\
    256$^2$ & 1.56e-4 &  2.00e-05  & 2.05  &  6.06e-05   & 1.23\\ \hline
  \end{tabular}
  \caption{Convergence of the temperature field, melting solid}
  \label{tab:planar1D_fixed}
\end{table}

\subsection{Frank's Spheres}
\begin{figure}[h!]
\begin{minipage}[h!]{.49\linewidth}
\includegraphics[width = 0.95\textwidth]{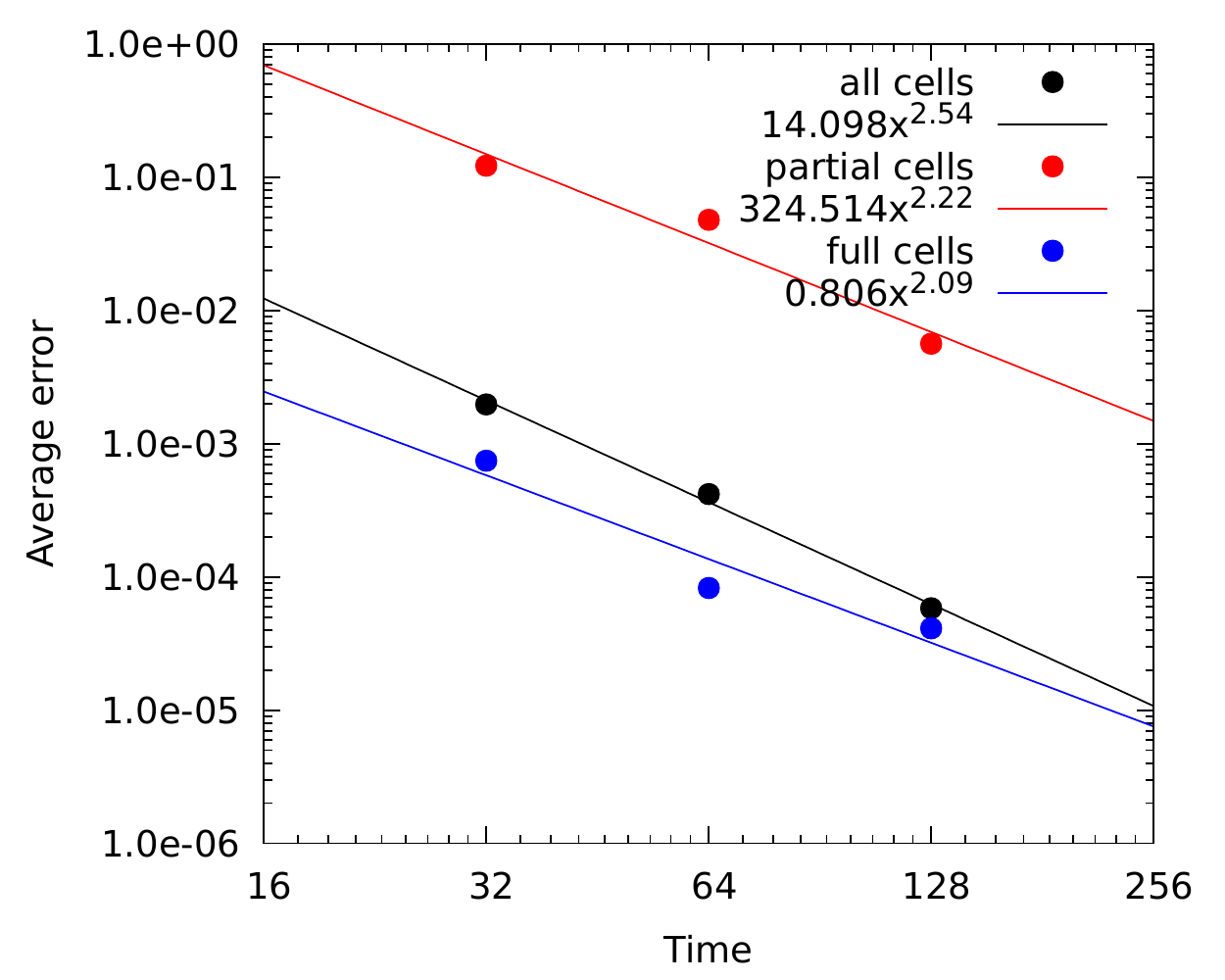}
  \subcaption{Frank's Spheres convergence results, fixed timestep $\Delta t = 10^{-4}$.} 
  \label{fig:FrankErrorFixed}
\end{minipage}
\begin{minipage}[h!]{.49\linewidth}
\centering
\includegraphics[width = 0.95\textwidth]{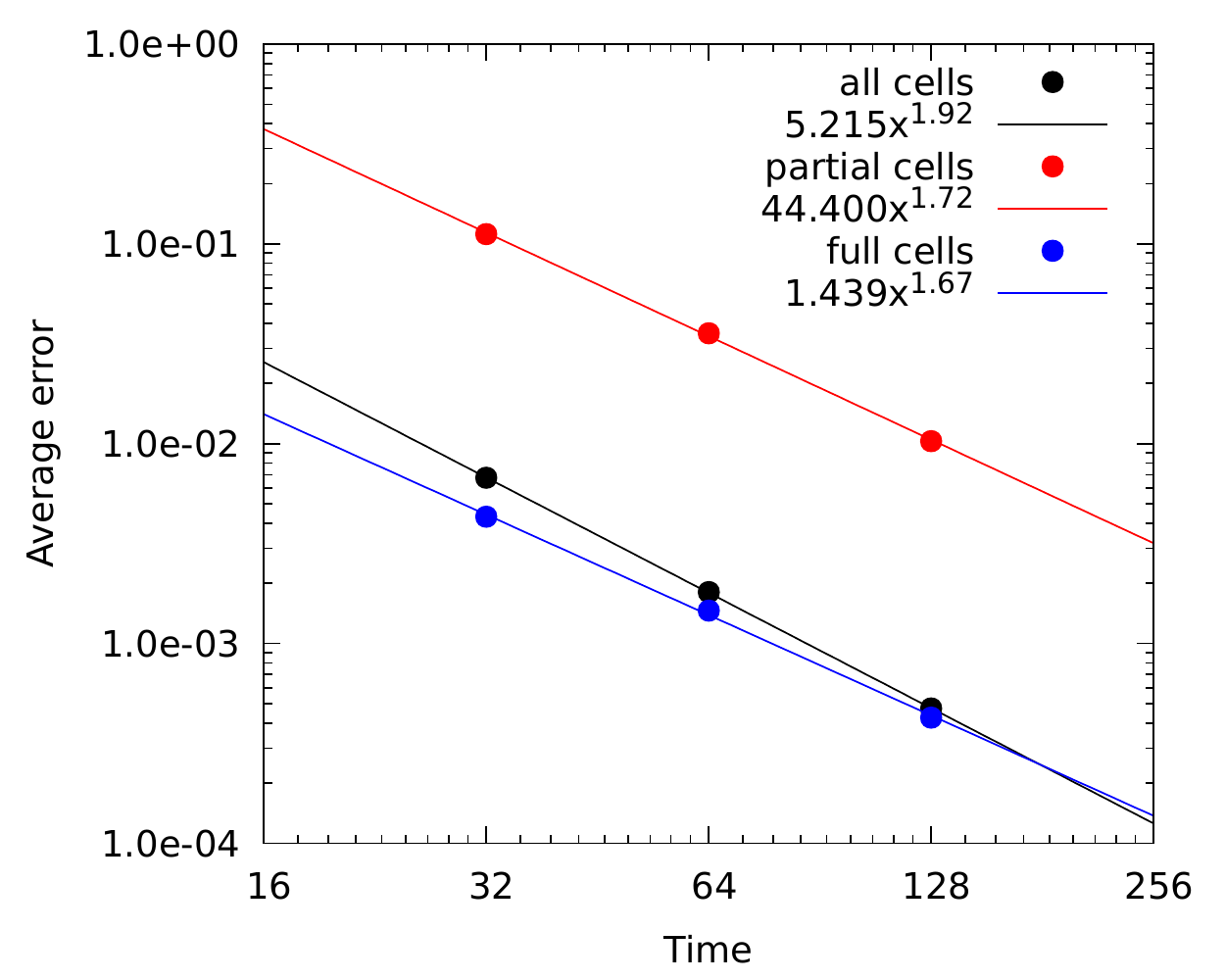}
  \subcaption{Frank's Spheres convergence results varying timestep and mesh size.}
  \label{fig:FrankError}
\end{minipage}
\centerline{\begin{minipage}[c]{.45\linewidth}
\centering
\includegraphics[width=0.9\textwidth]{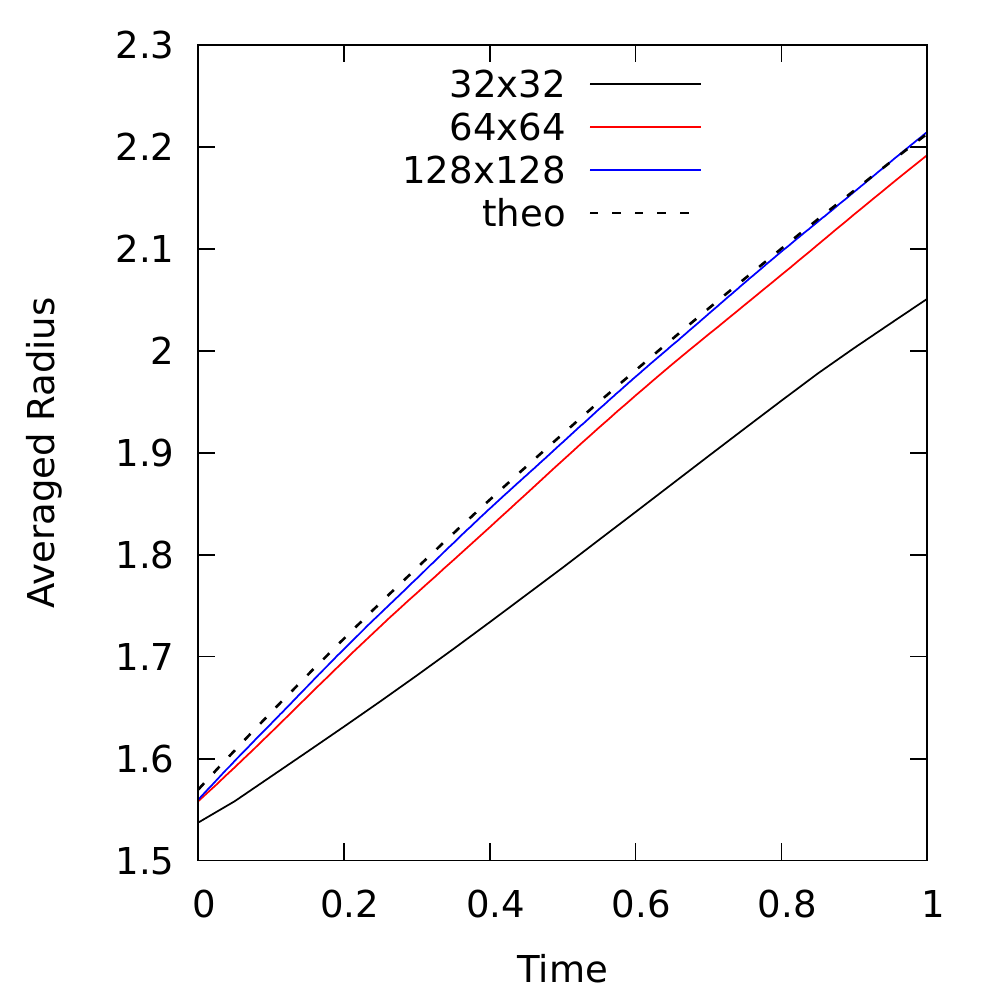}
\subcaption{Comparison of simulated radii with different grid size and
    theoretical prediction.}
\label{fig:FrankPosInterf}
\end{minipage}}
\caption{Results for Frank's spheres.}\label{fig:frank_global_results}
\end{figure}

Frank's spheres correspond to the growth of an ice sphere in an undercooled liquid.
The theory of this test case has been studied originally by Frank \cite{frank1950radially}, and it is a crucial test of the numerical stability of the scheme.
Indeed, in the absence of anisotropy, a growing sphere (whatever the 
space dimension $D=2$ or $3$ in practice) is an exact solution of the dynamics, although it is unstable due to the well known Mullins-Sekerka instability. The next cases in this paper focus on the simulation of dendritic growth. 
In the simulations, we want to stress that because of the numerical noise, the sphere destabilizes and forms dendrites: the numerical robustness of the scheme can be thus tested by investigating how the ice domain diverges from the sphere.
Therefore, starting with a spherical initial interface (circle in 2D or sphere
in 3D) containing a solid seed surrounded by an undercooled liquid, we will test the stability of the numerical scheme by inspecting the sphere growth. In fact, a class of
self-similar solutions has been developed by Frank, in one, two and three
dimensions, for which the sphere growth follows a square-root-of-time evolution characterized by the number $S$
\begin{equation}
R(t) = S t^{1/2}
\end{equation}
The solutions of the problem can be parametrized using the self similar variable $s=r/R(t)$.
and the corresponding dimensionless temperature field is $0$ for $s<S$  (corresponding to $r<R(t)$ and using for simplicity the same notation $T$ for the temperature field and its self-similar function):
\begin{equation}
T(r,t) = T(s) = T_\infty \left(1- \frac{F_D(s)}{F_D(S)}\right)
\end{equation}
if $s> S$. The functions $F_D$ are solutions of the equations and for $D=2$ we have:
\begin{equation}
F_2(s) = E_1\left(\frac{s^2}{4}\right)
\end{equation}
where
\begin{equation}
E_1(x) = \int_1^\infty \frac{e^{-xt}}{t}dt = \int_x^\infty \frac{e^{-t}}{t}dt
\end{equation}
Calculations are performed with the parameters of Almgren~\cite{Almgren1993}:
\begin{equation}
T_{\infty} = -0.5 = \frac12 S \dfrac{F_{2}(S)}{F'_{2}(S)}
\end{equation} 
which gives a value of $S=1.56$. Results of the calculation are shown on \cref{fig:frank_global_results} the initial time is $t_{0} = 1$. \cref{fig:FrankErrorFixed} shows results for calculation after 100 iterations with fixed timestep $\Delta t = 10^{-4}$, we recover a second-order convergence for the final temperature field. Calculations have also been performed for varying mesh size and timestep where $\Delta t = 0.2 (\Delta x)^2$. The order of accuracy is between $1.5$ and $2$.

\subsection{Crystal growth}
As discussed above, crystal formation is physically unstable and leads to
dendritic growth~\cite{Langer1980}. In order to study the formation of
dendrites we consider an ice crystal growing in an undercooled liquid as in
Chen \textit{et al.} \cite{chen_simple_1997}. In the absence of liquid flow,
this configuration consists in the diffusion of two tracers separated by a
complex embedded boundary. The dimensionless computational domain is $\Omega =
[-2:2]\times [-2:2]$, the initial interface (0-level-set) is defined by:
$\Gamma(r,\theta) = {(r,\theta) / r^2(1-0.3\cos (4\theta))-\dfrac{1}{15}}$
where $\theta$ is the angle of the outward normal with the $x$ axis. The
interface moves according to the Stefan relation with ${\rm St} = 0.5$. The
boundaries are thermally isolated, therefore we should reach a steady state
around the time $0.8$ with about half the domain that is solid.

\begin{figure}[h!]
\begin{minipage}[h!]{.33\linewidth}
\includegraphics[width= 0.9\textwidth]{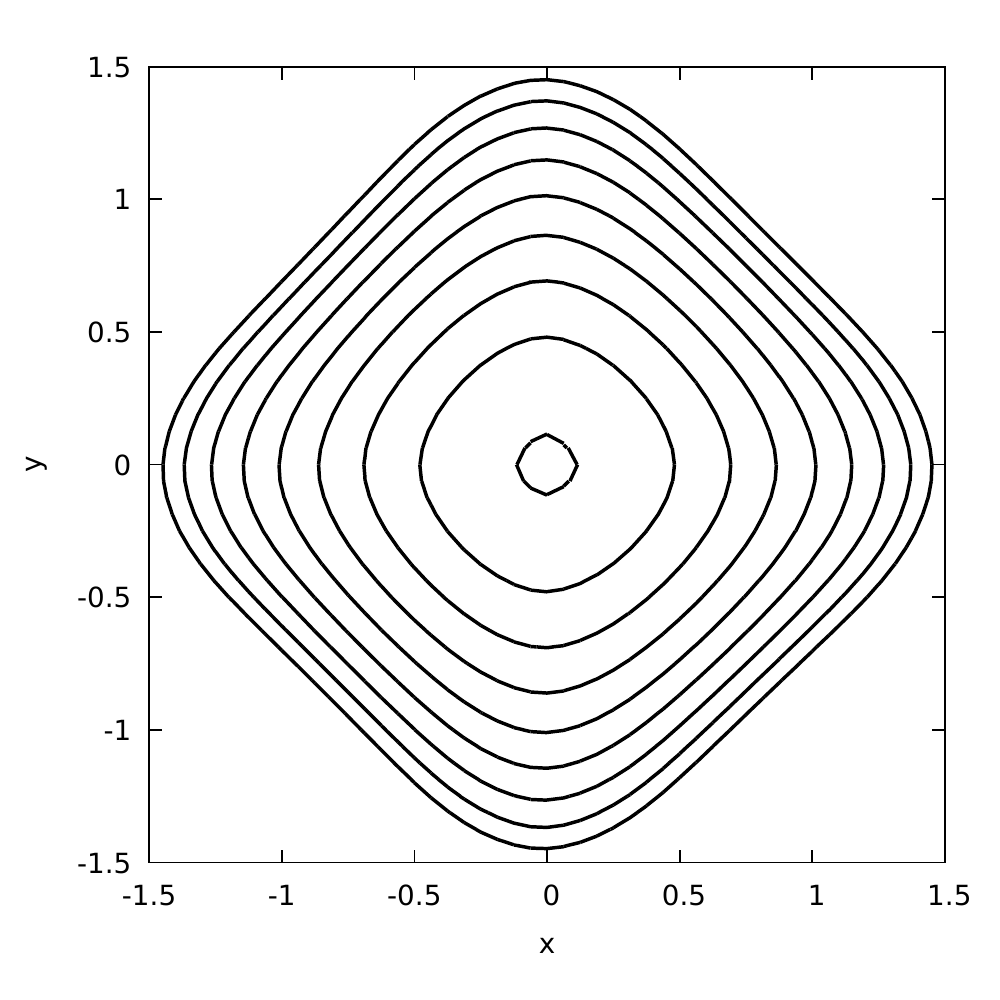}
\subcaption{$64^{2}$}
\end{minipage}
\begin{minipage}[h!]{.33\linewidth}
\includegraphics[width= 0.9\textwidth]{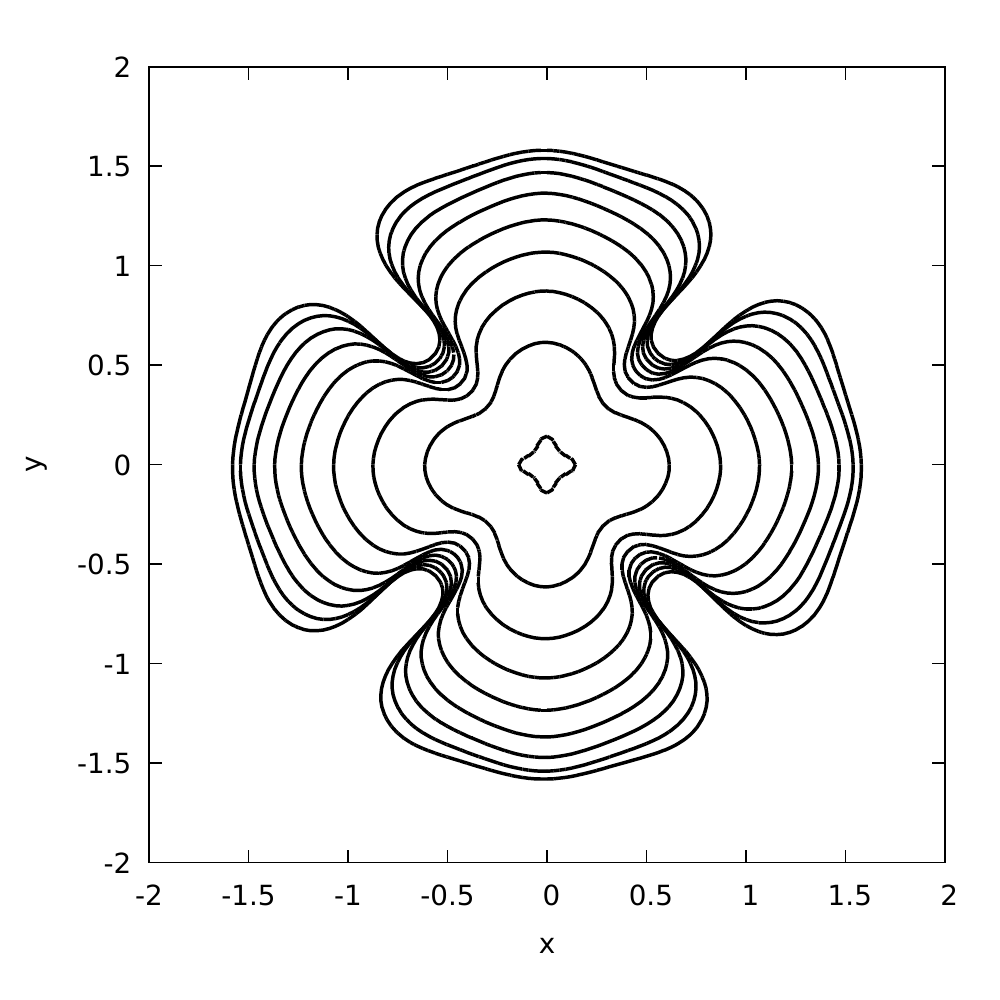}
\subcaption{$128^{2}$}
\end{minipage}
\begin{minipage}[h!]{.33\linewidth}
\includegraphics[width= 0.9\textwidth]{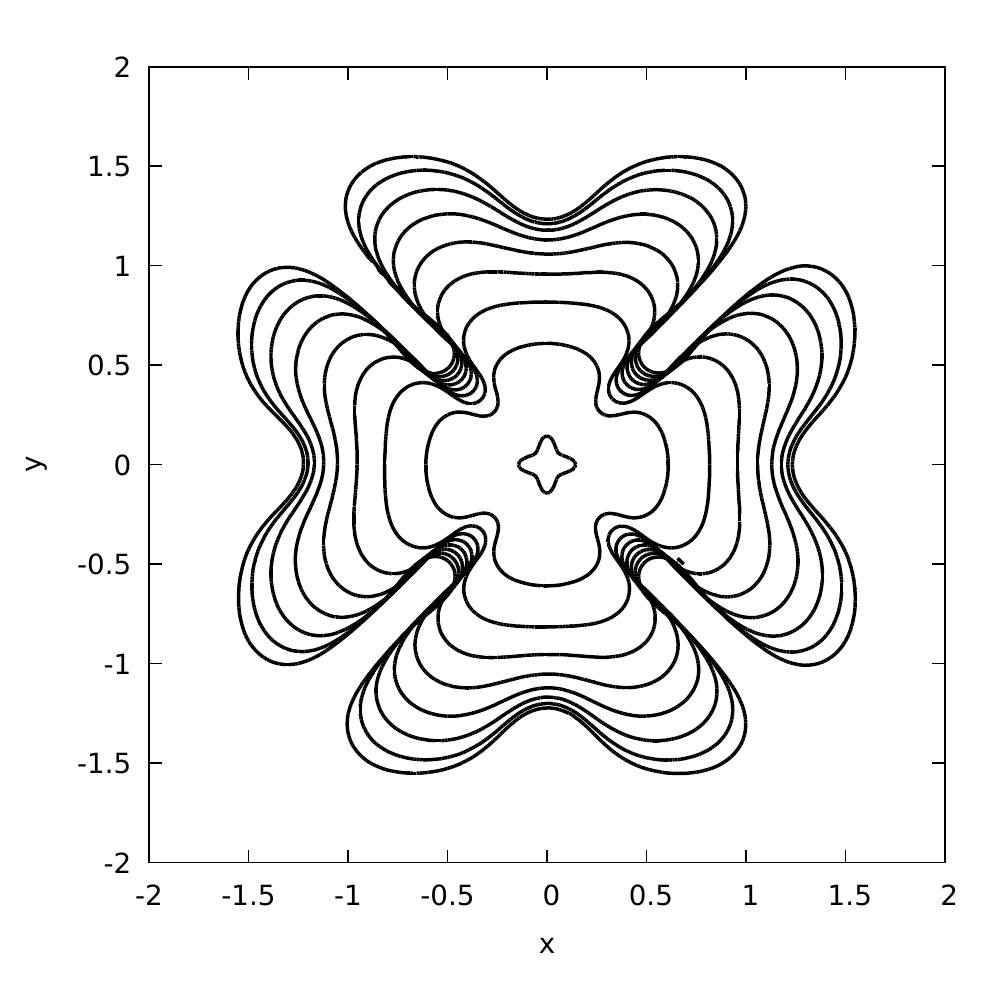}
\subcaption{$256^{2}$}
\end{minipage}
\caption{Influence of the spatial resolution on the initial dendritic growth. As the resolution increases, shorter wavelengths are resolved leading to an earlier development of the instability.}
\label{fig:crystal_growth}
\end{figure}

The ice particle is initially at $T_S = 0$ and the temperature in the liquid
is $T_L = -0.5$. The temperature on the interface follows the
Gibbs--Thomson relation, taking $\epsilon_{\kappa} = \epsilon_{v} = 2\times
10^{-3}$. The interface is plotted every 0.1 unit time unit on 
\cref{fig:crystal_growth} for three different grid sizes, showing an instability which depends on the grid size: in fact the instability generates high-curvature unstable regions that are eventually stabilized by the Gibbs--Thomson contribution in the melting temperature. Therefore, the smaller the mesh size, the faster and the more complex the instability grows. The length of the dendrites can be directly calculated and is fixed by the value of $\epsilon_{\kappa}$.
Our results are in fact quite comparable with \cite{chen_simple_1997,Tan2006}, but the onset of instabilities can be seen on the $256^{2}$ case far earlier than in their simulations, indicating a very low level of built-in regularization in our method.

\subsection{Crystal growth with sixfold anisotropy}\label{sec:sixfold}
It is well known that to describe accurately dendritic growth in crystals, anisotropy of the Gibbs--Thomson condition on the curvature needs to be implemented, taking for instance the form~\cite{Tan2006}, in 2D:
\begin{equation}
\epsilon_{\kappa} = 0.001\left(1+\epsilon\left[\dfrac{8}{3}\sin^4
(3(\theta-\frac\pi2))-1\right]\right).
\end{equation}
We will use $\epsilon = 0.4$, $\theta = (Ox,\boldsymbol{n})$ in our simulations. We expect to have 6 primary dendrites growing at the same speed, details of the calculation are given in \cref{tab:sixfold}.

\begin{table}[h!]
  \centering
  \begin{tabular}{|c|c|c|}
  \hline
  Undercooling & Domain size & $\epsilon_{v}$ \\
  \hline
  0.8  & $[-2:2]\times[-2:2]$ & 0.001 \\
  \hline
  \end{tabular}
  \caption{Details of the calculation}\label{tab:sixfold}
\end{table}

\begin{figure}[h!]
  \centering
  \includegraphics[width= 0.5\textwidth]{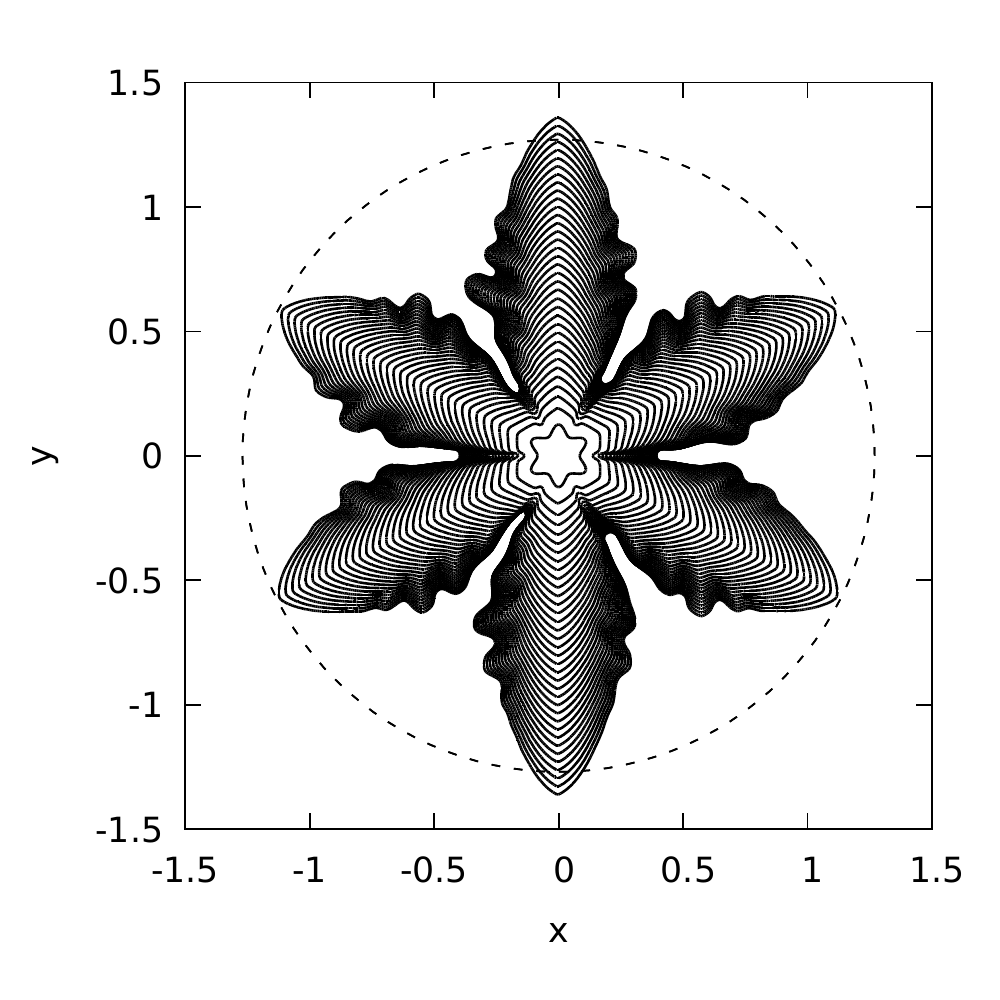}
  \caption{Anisotropy effect with a sixfold symmetry}
  \label{fig:sixfold}
\end{figure}

We show results of the calculations on \cref{fig:sixfold}, where the interface is plotted every $\Delta t = 3\times10^{-3}$ and the final time of the calculation is $t = 3.6\times10^{-2}$. Since we use adaptive mesh refinement, the maximum equivalent resolution is $512^{2}$. We plot a circle of radius 1.27 with a dashed line, which is a simple fit to compare the size of the dendrites. At the final time, the size of the six main dendrites is very similar which is a strong indication of the correct treatment, even in non-grid-aligned directions, of the anisotropic Gibbs--Thomson condition. We notice also the presence of secondary dendrites, as expected from crystal growth instability.

\subsection{Crystal growth in 3D}

\begin{figure}[h!]
\begin{center}
\begin{minipage}[h!]{.4\linewidth}
\centering
\includegraphics[width=0.99\textwidth]{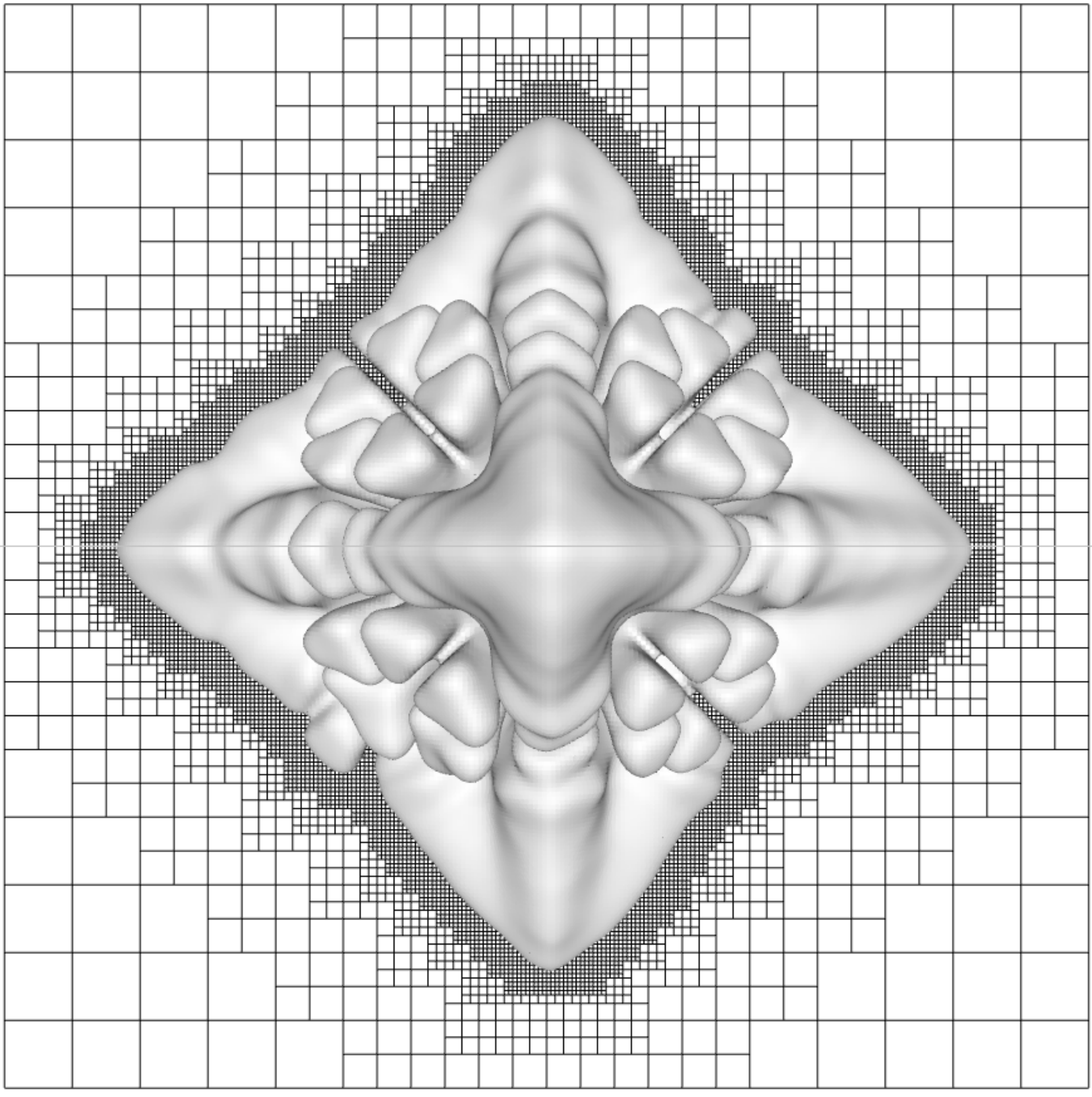}
\end{minipage}
\begin{minipage}[h!]{.4\linewidth}
\centering
\includegraphics[width=0.99\textwidth]{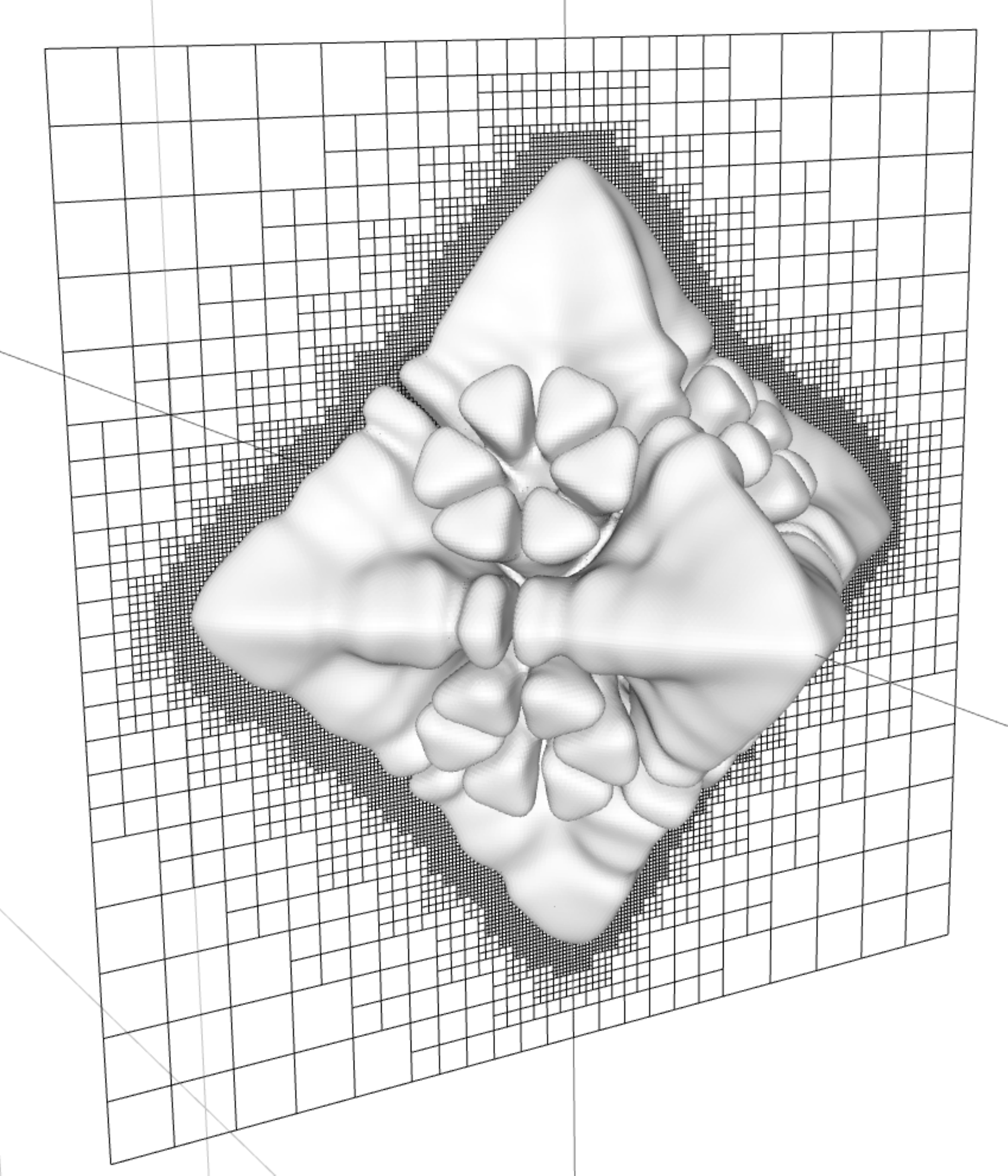}
\end{minipage}
\caption{Crystal growth in 3D}
\label{fig:Crystal3D}
\end{center}
\end{figure}

The details of the calculation are similar to those of the previous calculation, using now a 3D Gibbs--Thomson condition, following:
\begin{equation}
   \epsilon_{\kappa} = \overline{\epsilon_{\kappa}}\left[1 -3\epsilon_{4}  + 4\epsilon_{4}^{2}\left(\sum_{i=1}^{3}n_{i}^{4}\right)\right]
\end{equation}
with $\epsilon_{4} = 0.4$, $\overline{\epsilon_{\kappa}} = 0.001$, the initial undercooling is $0.8$ as in the previous calculation. Here the results are shown for a grid with an equivalent resolution of $256^{3}$, the final number of cells is about $2\times 10^{6}$ and the final time is $3.6\times10^{-2}$. The value $\epsilon_{4}=0.4$ introduces a strong anisotropy on the interface temperature. \cref{fig:Crystal3D} represents the interface at the end of the calculation and a slice of the mesh in a medial plane. Results are quite similar to the simulations by Lin \textit{et al.} \cite{lin2011adaptive}. One can see the expected fourfold periodicity of the main dendrites. The secondary dendrites are also quite well captured by our method demonstrating both its robustness regarding anisotropy and its accuracy. We point out the strength of the mesh adaptation method allowing very local mesh refinement. Future simulations will focus on having a locally converged state with regards to the mesh adaptation criteria. The main driving adaptation criterion is linked to the thermal boundary layer formed near the interface.

\subsection{Comparison of the tip velocity with linear solvability theory}

\begin{figure}[h!]
\centerline{\begin{minipage}[c]{0.6\linewidth}
  \includegraphics[width=\linewidth]{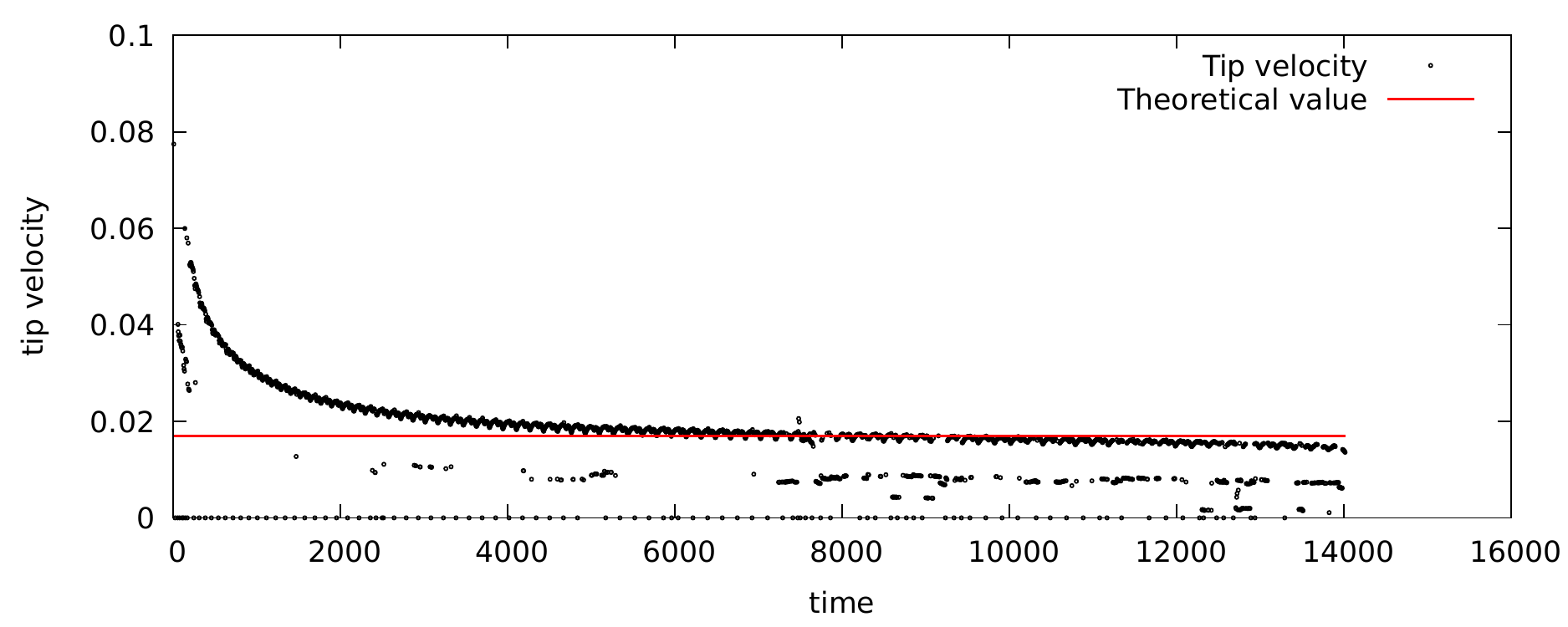}
  \subcaption{Whole simiulation}
\end{minipage}}
\centerline{\begin{minipage}{0.6\linewidth}
  \includegraphics[width=\linewidth]{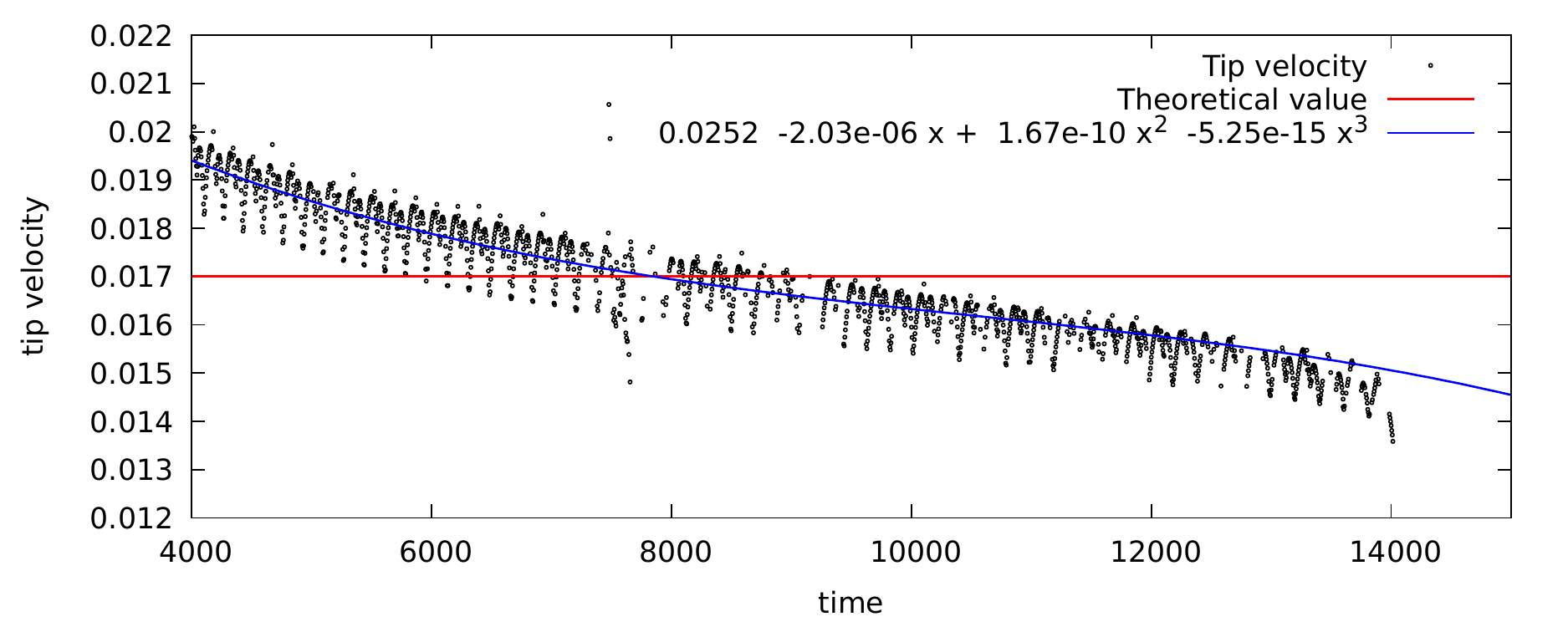}
  \subcaption{Zoom. Oscillations in the interface velocity}\label{subfig:linear_solv_oscillations}
\end{minipage}}
\caption{Tip velocity as a function of time}\label{fig:speed_linear_solvability}
\end{figure}

This case models anisotropy effects and can be compared with predictions of microscopic solvability theory \cite{Tan2006,kim2000computation}. Here the Gibbs--Thomson condition is:
\begin{equation}
T_{\Gamma} = -\overline{\epsilon_{\kappa}}(1 - 15 \epsilon \cos 4 \theta)\kappa
\end{equation}
with $\overline{\epsilon_{\kappa}} = 0.5$, $\epsilon = 0.05$, the initial
undercooling is $0.55$ and every other parameter is unity. The mesh
has an equivalent resolution of $512^{2}$ and the computational domain is $
[-400:400]\times[-400:400]$. The results are plotted for the
adimensionalized field $\tilde{V}$, $\tilde{x}$, and $\tilde{t}$ see \cref{fig:speed_linear_solvability}. A difficulty associated with this test case is the initialization of the temperature in the newly solid cells which depends on the temperature of the interface and therefore its velocity with the Gibbs--Thomson condition. As expected, small oscillations in the temperature gradients influence the motion of the interface and create oscillations in the velocity of the interface, \cref{subfig:linear_solv_oscillations}. The expected tip velocity is $1.7\times{10^{-3}}$, in our simulation, the velocity reaches the value $1.9\times10^{-3}$ then drifts slowly, as shown by the regression coefficients in the \cref{subfig:linear_solv_oscillations}. The global drift is probably due to the influence of the boundary. The local oscillations are linked to i)the jumps in discretization stencil and the associated truncation errors and ii)the fact that we solve fixed-mesh problems at each timestep, something that could be fixed by reformulating our method using a Arbitrary Lagrangian Eulerian framework.

\begin{figure}[h!]
\begin{minipage}[]{0.6\linewidth}
  \includegraphics[width=\linewidth]{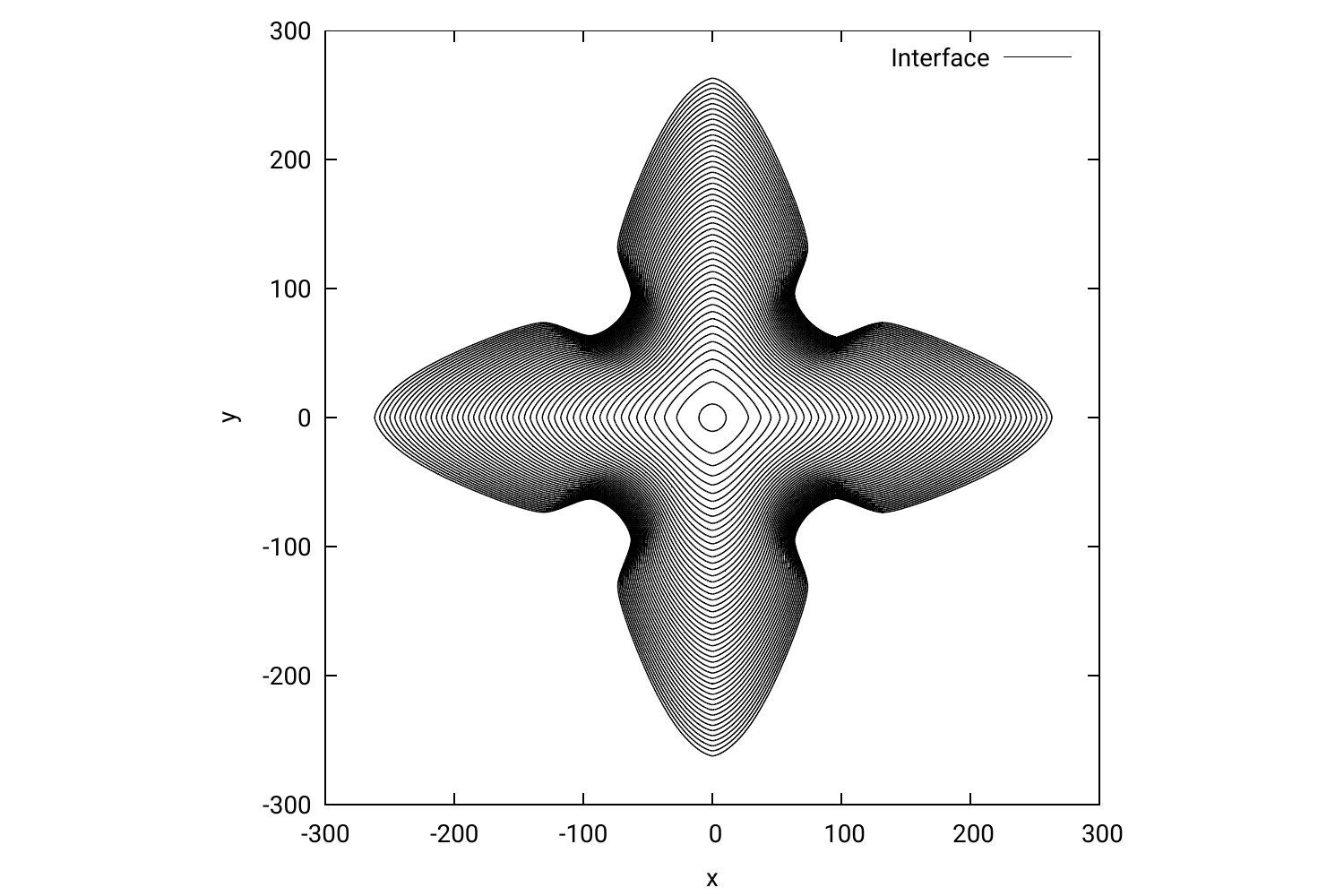}
  \subcaption{Interface output at different instants during the simulation}
\end{minipage}
\begin{minipage}[]{0.5\linewidth}
  \centering\includegraphics[width=0.75\linewidth]{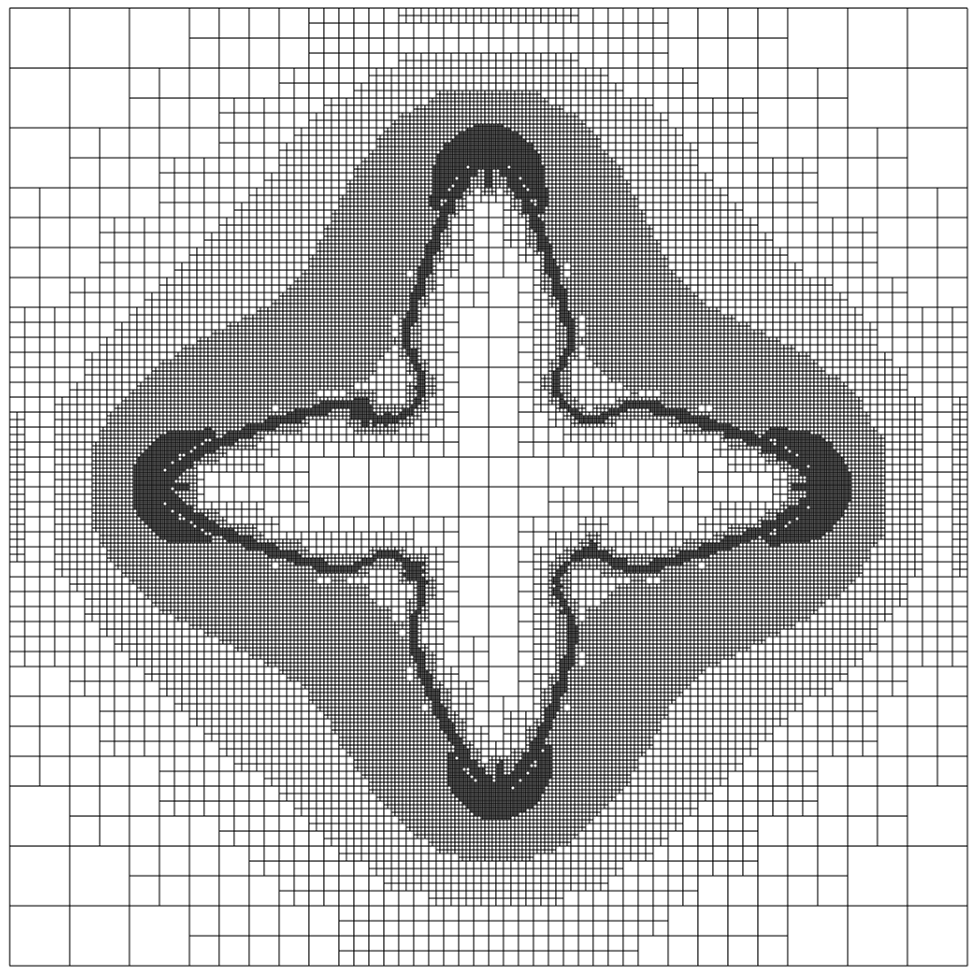}
  \subcaption{Mesh at final time}
\end{minipage}
  \caption{Linear solvability test case, interface during calculation}
\end{figure}

\subsection{Rayleigh--B\'enard instability with a moving melting boundary}
\begin{figure}[h!]
\centering
\includegraphics[width=0.5\textwidth]{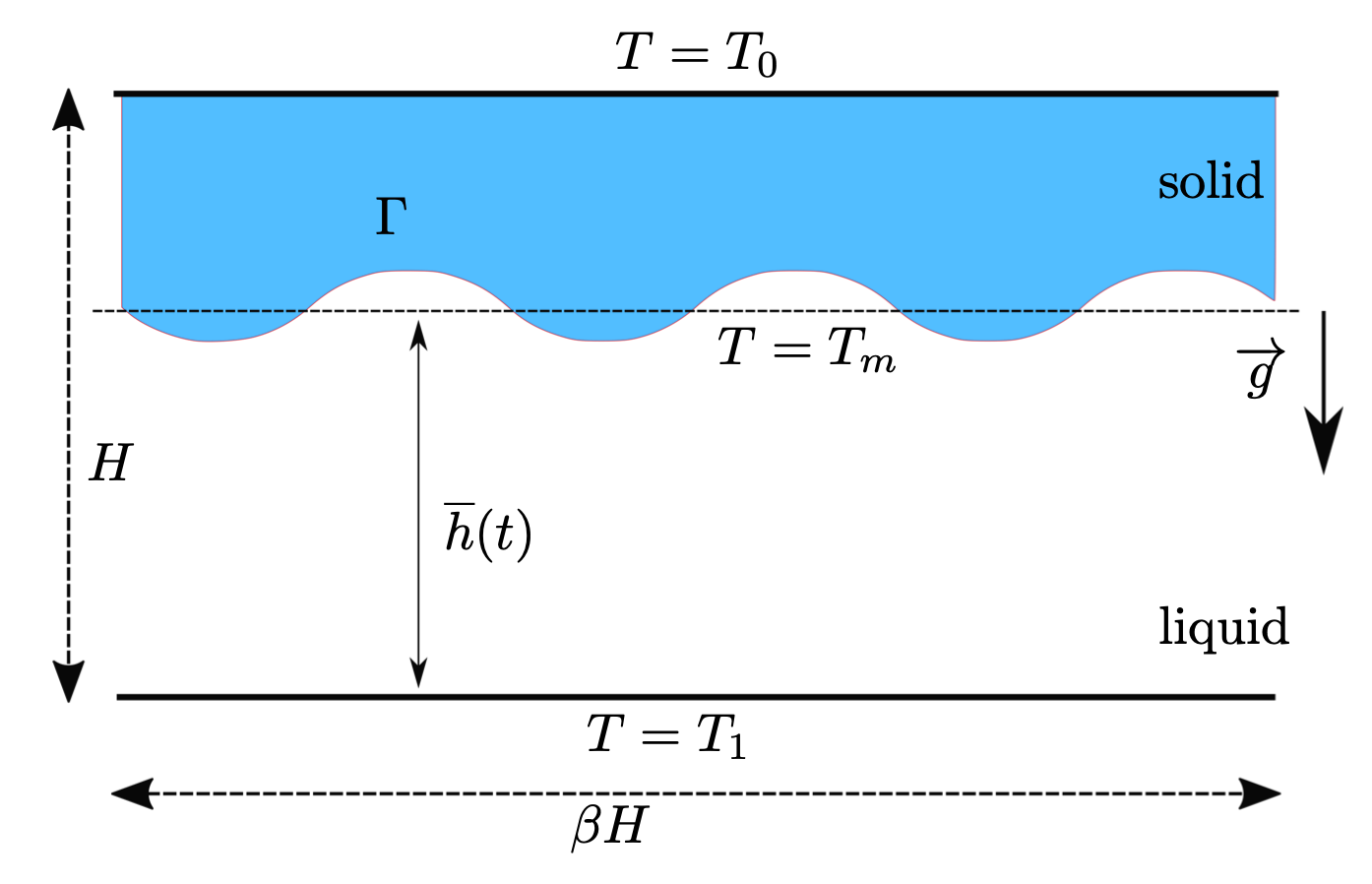}
      \caption{Scheme of the Rayleight-Benard instability with a melting
      boundary}\label{fig:scheme_RB}
\end{figure}

Finally, in order to validate the coupling of our new method for solving
solidification fronts with the Navier-Stokes equation, we study the threshold
of the Rayleigh--B\'enard instability of a melting ice layer, following the
recent study of Favier \textit{et al.} \cite{Favier2019}. The configuration
studied is depicted on \cref{fig:scheme_RB}. A pure and incompressible
material under the influence of gravity $\boldsymbol{g} = -g
\boldsymbol{e_{Z}}$ is comprised between two walls such that it is heated from
below (by imposing a temperature $T_1$ at $z=0$) and cooled from above ($T_0$
at $z=H$). The melting temperature varies between these two imposed
temperatures $T_0<T_m<T_1$, so that, taking $\lambda_L=\lambda_S$, the
equilibrium position of the ice layer is simply determined by the balance of
the thermal fluxes, giving: $$h_m=H\frac{T_m-T_0}{T_1-T_0}.$$ This
configuration is similar to the classical one for Rayleigh--Benard
(R-B)instabilities, but this time the upper boundary of the liquid can move
through melting or freezing and thus initiate the R-B instability during the
dynamics. Depending on the physical parameters and the initial conditions,
different stationary regimes have been observed in numerical simulations using
phase-field modeling for phase change\cite{purseed2020bistability}.

We thus consider a two-dimensional system bounded by two horizontal walls as shown in \cref{fig:scheme_RB} separated by a distance $H$ with periodic boundary conditions on the horizontal direction with an aspect ratio $\beta$ between the horizontal and vertical dimensions. We apply no slip boundary conditions on the upper and lower boundaries and also on the interface between the two phases $\Gamma$. We perform the simulation using the Boussinesq approximation, where the variation of the liquid density with the temperature is taken into account in the buoyancy force only. More precisely, we assume that the density $\rho_S=\rho_L=\rho$ and the thermal diffusivity $D_{S}(T)=D_L(T)=D_0$ are constant and equal in both domains, as well as the fluid viscosity. The dimensionless set of equations thus reads in the fluid domain:
\begin{eqnarray}
\frac{1}{\sigma} (\frac{\partial \boldsymbol{u}}{\partial t} + \boldsymbol{u}\cdot \nabla 
\boldsymbol{u}) &=& - \nabla P + Ra\,\theta\,\boldsymbol{e_{y}} + \nabla^{2}\boldsymbol{u}\\
\frac{\partial \theta}{\partial t} + \boldsymbol{u}\cdot \nabla \theta&=&\nabla^{2}\theta
\end{eqnarray}
where $\theta = \dfrac{T-T_{0}}{T_{1}-T_{0}}$ is the dimensionless reduced temperature, $P$ is the dimensionless pressure, $Ra$ and $\sigma$ are the Rayleigh and Prandtl numbers respectively:
\begin{equation}
Ra = \dfrac{g\alpha(T_{1}-T_{0})H^{3}}{\nu D_{0}} \text{ and }\sigma = \dfrac{\nu}{D_{0}} \text{ .}
\end{equation}
where $\alpha =1$ is the thermal expansion coefficient. As in \cite{Favier2019}, we impose $\sigma = 1$ and only the Rayleigh number is varied for all of our simulations. In the solid phase we have:
\begin{equation}
\dfrac{\partial \theta}{\partial t} = \nabla^{2} \theta \text{ .}
\end{equation}
The Stefan condition is applied on the boundary, and the temperature on the interface is supposed to be constant $\theta_{\Gamma} = \theta_{m}$. We apply our method with two different solvers, one for the Navier--Stokes equations in the liquid and a simple diffusion solver in the solid.

\begin{figure}[h!]
\centering
\includegraphics[width=0.6\textwidth]{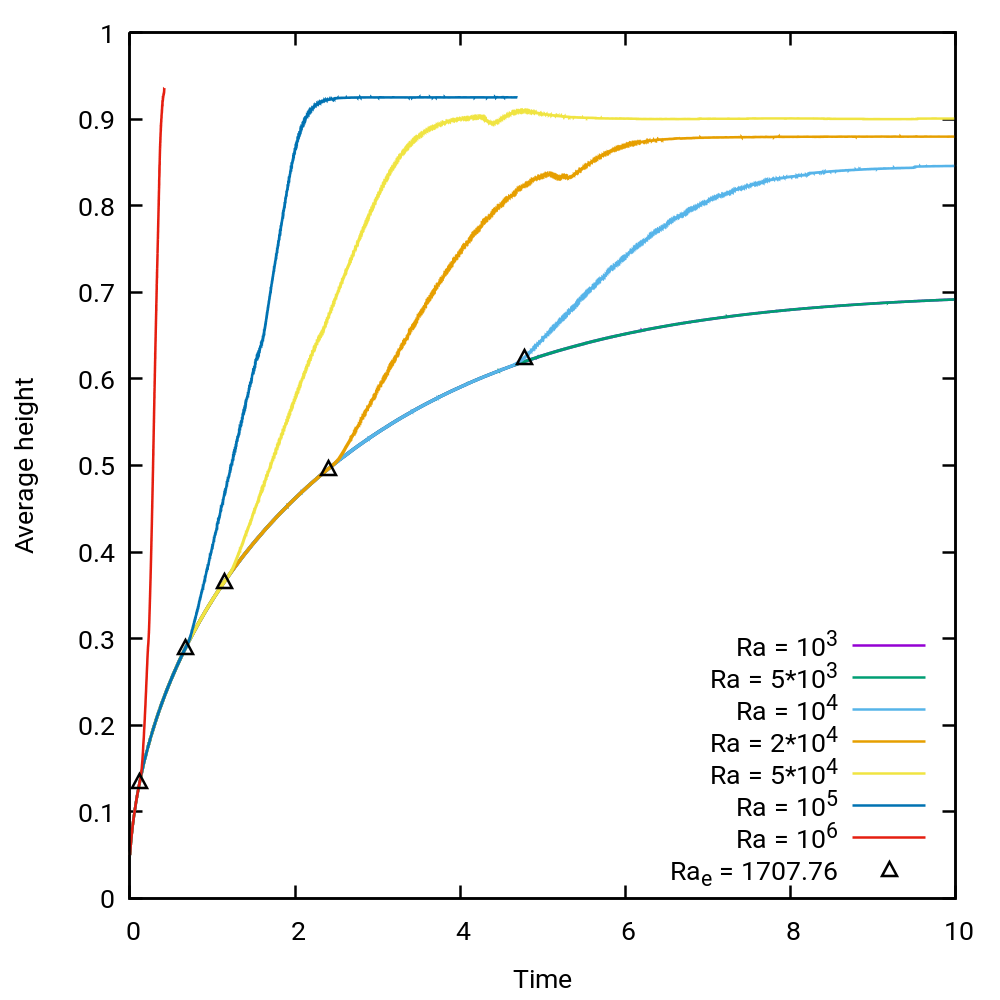}
\caption{Rayleigh--B\'enard instability with a melting boundary. Average height evolution for different Rayleigh numbers}\label{subfig:heigh_RB}
\end{figure}

We re-perform the calculations of Favier to study the onset of the R-B instability and the formation of convection cells. We define similarly the effective Rayleigh number defined using the fluid layer thickness following:
\begin{equation}
Ra_{e} = Ra (1-\theta_{m}) \overline{h}^{3}
\end{equation}
where the averaged fluid height $h(t)$ is defined as
\begin{equation}
\overline{h}(t)=\dfrac1\beta \int_{0}^{\beta}h(x,t)dx.
\end{equation}
Convection cells are expected to appear once the simulations reach a critical effective Rayleigh number $Ra_{c} = 1707.76$. Details of the initial grids are given in \cref{tab:tab1}: the initial effective Rayleigh number $Ra_{e}$ is much lower than the critical Rayleigh number such that the calculations always start with a diffusion-driven dynamic. Results of the average height during the calculations are given in \cref{subfig:heigh_RB}. For each curve, we added a triangle sign, to position when the simulation reaches the critical Rayleigh number. This collection of triangles clearly separates two regimes, the diffusion-driven phase from the convection-influenced one. Once the apparent critical Rayleigh number $1707.76$ is reached, the thermal exchange between the bottom boundary and the interface is greatly enhanced and the interface melts much faster.

\begin{table}[h!]
\centering
\begin{tabular}{|c|c|c|c|c|c|l|}
\hline
$N_{x}$ & $N_{y}$ & $Ra$ & $St$ & $\theta_{m}$ & $\beta$ & $h_{0}$ \\ \hline
$512$ & $64$ & $[10^3:10^6]$ & $10$ & $0.3$ & $8$ & $0.05$ \\ \hline
\end{tabular}
\caption{Numerical set up for the study of the critical Rayleigh number $Ra_{c}$}
\label{tab:tab1}
\end{table}

\begin{figure}[h!]
\centering
\includegraphics[width=0.75\textwidth]{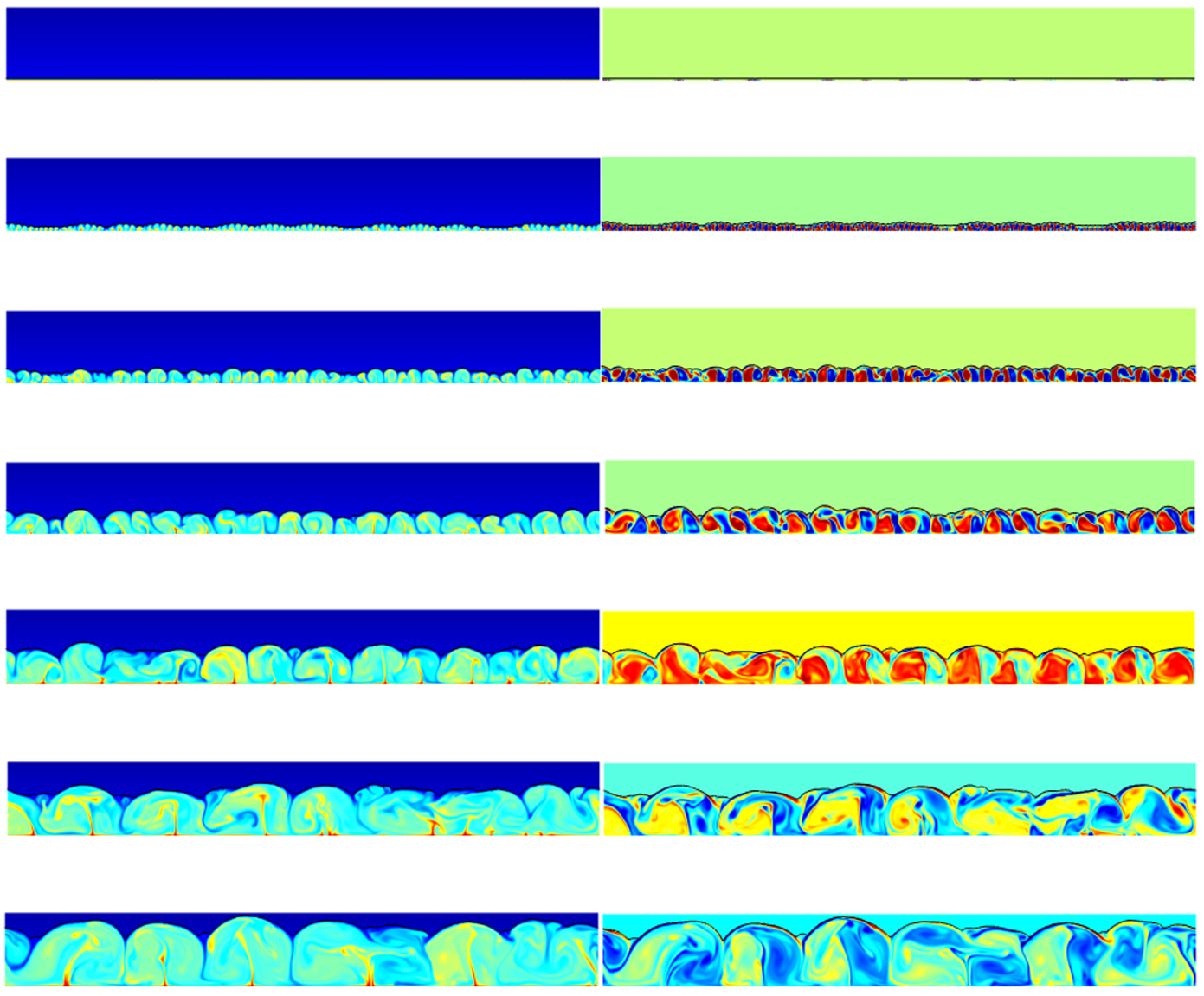}
\caption{$Ra = 1\times 10^{8}$ case, at instant: $t = 0.0005$, $t = 0.0011$, $t = 0.0016$, $t = 0.0023$, $t = 0.0033$, $t = 0.0044$ and $t = 0.0059$. Colored by left: temperature, right: vorticity.}\label{subfig:highRB}
\end{figure}

For sufficiently high Rayleigh numbers, the characteristic size of the convection cells of this flow will vary with a secondary bifurcation mechanism. Secondary bifurcations also occured for $Ra= 10^{5}$ and $10^{6}$ in the previous calculations. These bifurcations occur once the averaged height $\overline{h}$ is equal to the characteristic wavelength of the convection rolls. In the cases where the Rayleigh number is about $10^{5}$--$10^{6}$ once the secondary bifurcation is reached, these convection cells have a sufficient time to merge and re-stabilize because the motion of the melting boundary is sufficiently slow. We also performed a simulation with a higher Rayleigh number ${\rm Ra} = 10^{8}$, details of the calculation are given in \cref{tab:high_ra}. In the ${\rm Ra} = 10^{8}$ case, convection cells never fully stabilize giving birth to many unsteady thermal plumes as shown in \cref{subfig:highRB}.

\begin{table}[h]
\centering
\begin{tabular}{|l|l|l|l|l|l|}
\hline
 $N_{x}$ & $Ra$   & $St$ & $\theta_{M}$ & $\alpha$ & $h_{0}$ \\ \hline
$2048$   & $10^8$ & $1$  & $0.05$       &  $8$     & $0.02$  \\ \hline
\end{tabular}
\caption{Parameters for the high Rayleigh configuration}
\label{tab:high_ra}
\end{table}

\section{Conclusion}

An original level-set embedded boundary hybrid method has been developed for
the simulation of liquid-solid phase change. Its key features are (a) the use
of finite-volume conservative operators for embedding the interface which is
seen as a boundary from each phase's perspective, (b)the associated
second-order accuracy on the gradients on the boundary, (c) a simple velocity
extension method, (d) an initialization method derived from the embedded
boundary method. The method has been validated on numerous classical
melting/solidification problems, from planar melting to dendritic growth and
an extension of Rayleigh-B\'enard problem in the presence of phase change. The
method exhibits spatial and time convergence orders ranging between 1.5 and 2.
It has been validated in two and three space dimensions and the use of
adaptive mesh refinement allows a precise account of dendritic growth for
instance. Future work should develop the coupling of this method with a Volume
of Fluid one to allow three phases (gas-liquid-solid) simulations, with a
special attention to the contact line dynamics, but also the implementation of
the density variation between the liquid and solid phases.

\appendix

\section{Redistancing method of Min \& Gibou}
\label{sec:Redistancing}
\begin{figure}[h!]
    \begin{minipage}[h!]{0.33\textwidth}
  {\includegraphics[width=.9\textwidth]{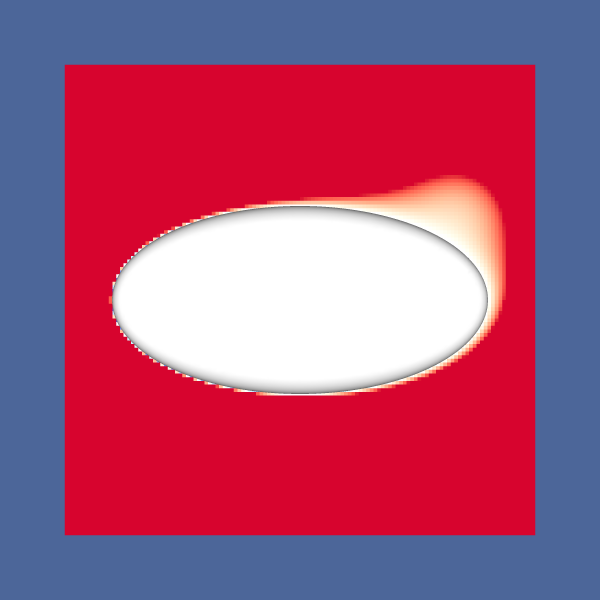}}
  \end{minipage}
  \begin{minipage}[h!]{0.33\textwidth}
  {\includegraphics[width=.9\textwidth]{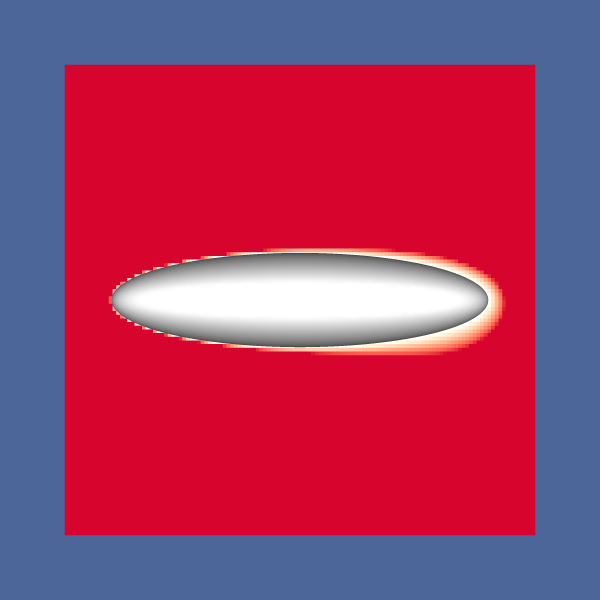}}
  \end{minipage}
  \begin{minipage}[h!]{0.33\textwidth}
  {\includegraphics[width=.9\textwidth]{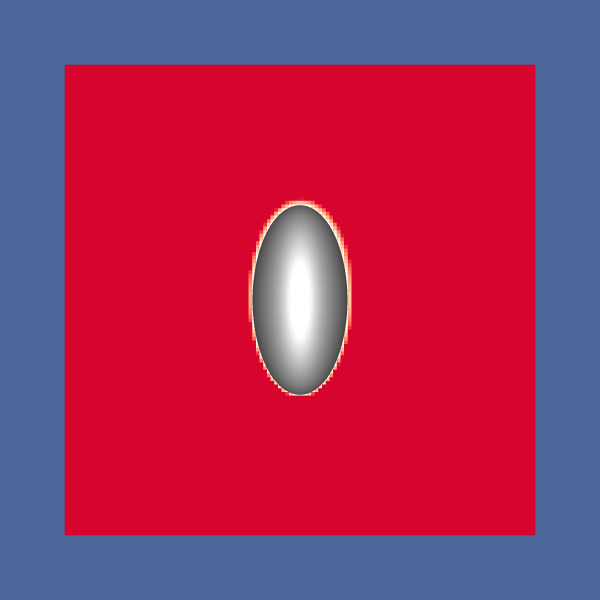}}
  \end{minipage}
  \caption{Initial value of the 0-level-set}
  \label{fig:reinit_0LS}
\end{figure}

\begin{figure}[h!]
    \begin{minipage}[h!]{0.33\textwidth}
  {\includegraphics[width=.9\textwidth]{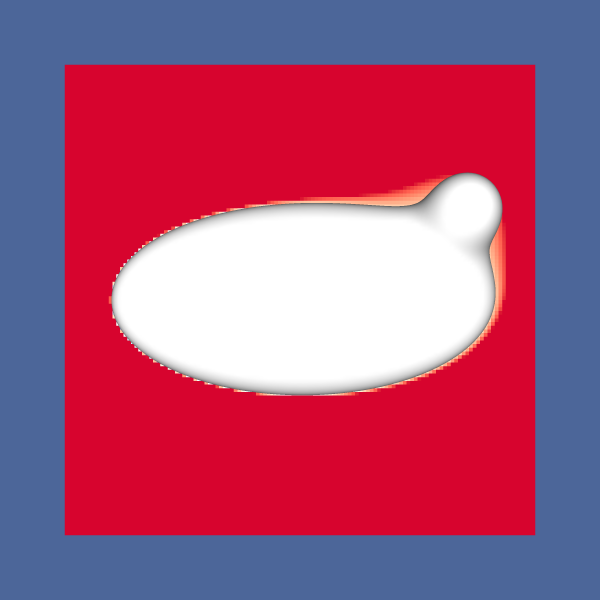}}
  \end{minipage}
  \begin{minipage}[h!]{0.33\textwidth}
  {\includegraphics[width=.9\textwidth]{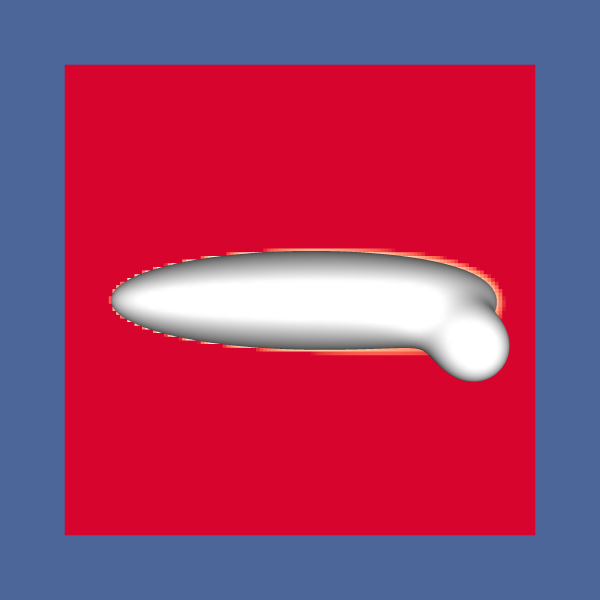}}
  \end{minipage}
  \begin{minipage}[h!]{0.33\textwidth}
  {\includegraphics[width=.9\textwidth]{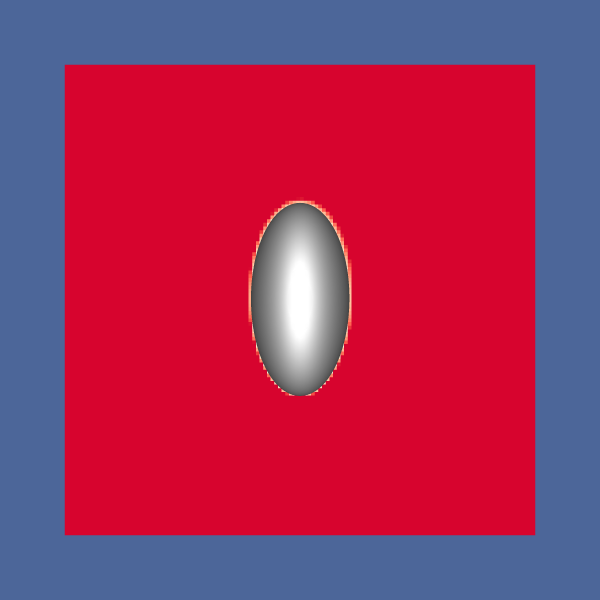}}
  \end{minipage}
  \caption{Initial value of the 0.8-level-set}
  \label{fig:reinit_0.8LS}
\end{figure}

\begin{figure}[h!]
    \begin{minipage}[h!]{0.33\textwidth}
  {\includegraphics[width=.9\textwidth]{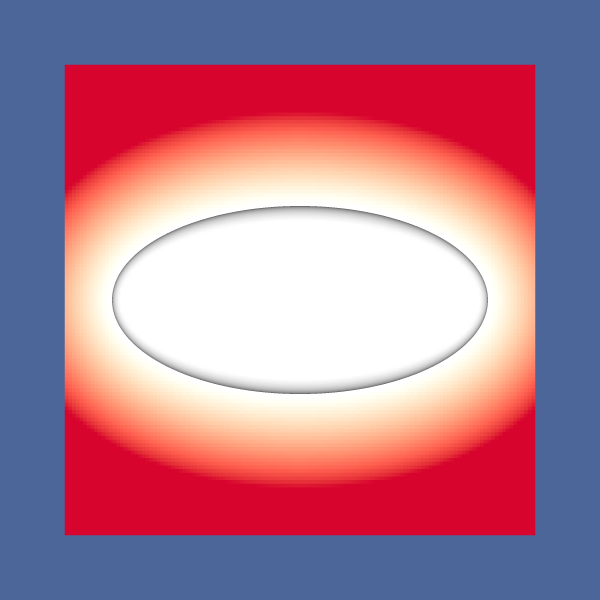}}
  \end{minipage}
  \begin{minipage}[h!]{0.33\textwidth}
  {\includegraphics[width=.9\textwidth]{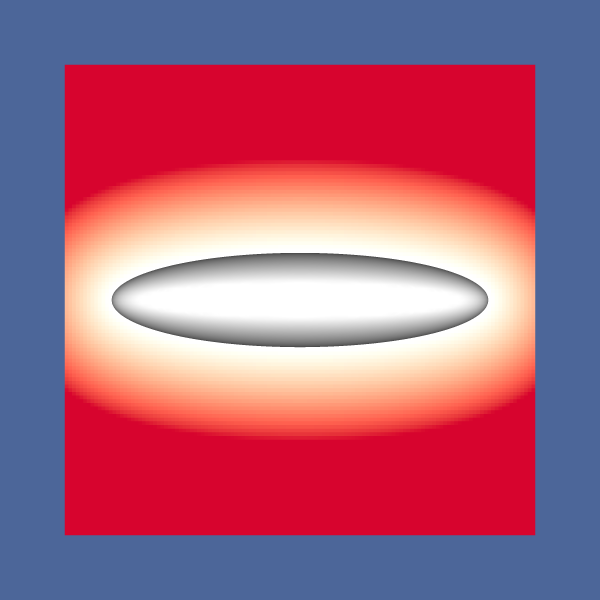}}
  \end{minipage}
  \begin{minipage}[h!]{0.33\textwidth}
  {\includegraphics[width=.9\textwidth]{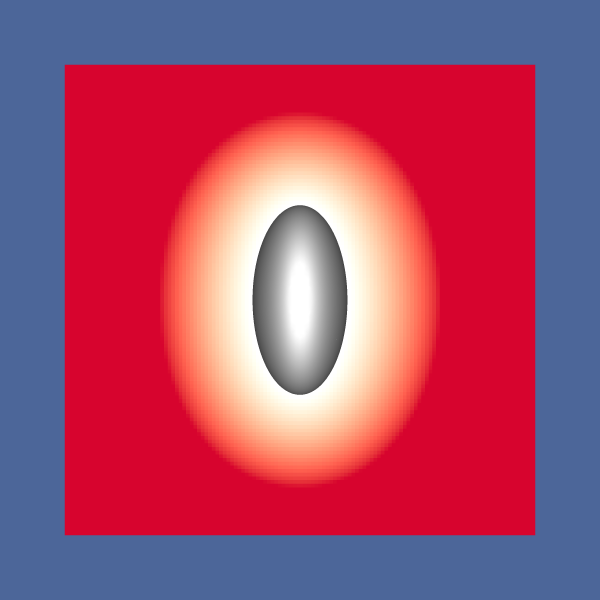}}
  \end{minipage}
  \caption{Final value of the 0-level-set}
  \label{fig:reinit_0LS_final}
\end{figure}

\begin{figure}[h!]
    \begin{minipage}[h!]{0.33\textwidth}
  {\includegraphics[width=.9\textwidth]{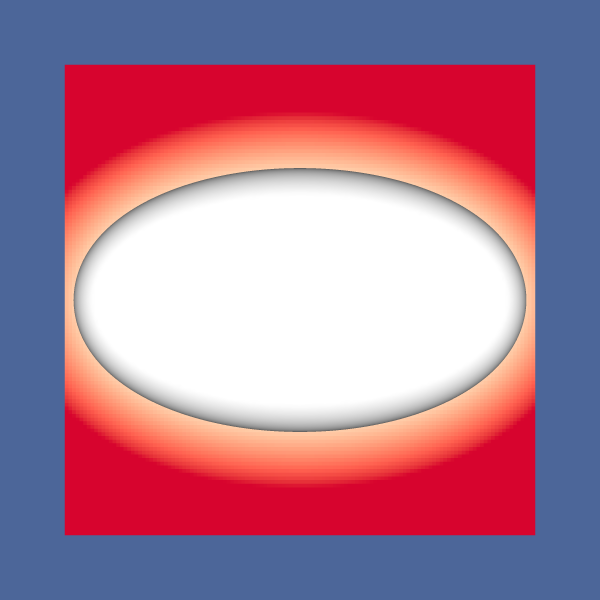}}
  \end{minipage}
  \begin{minipage}[h!]{0.33\textwidth}
  {\includegraphics[width=.9\textwidth]{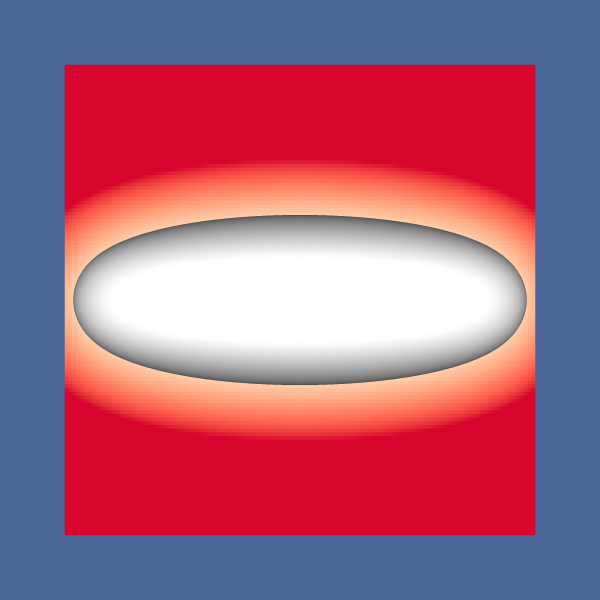}}
  \end{minipage}
  \begin{minipage}[h!]{0.33\textwidth}
  {\includegraphics[width=.9\textwidth]{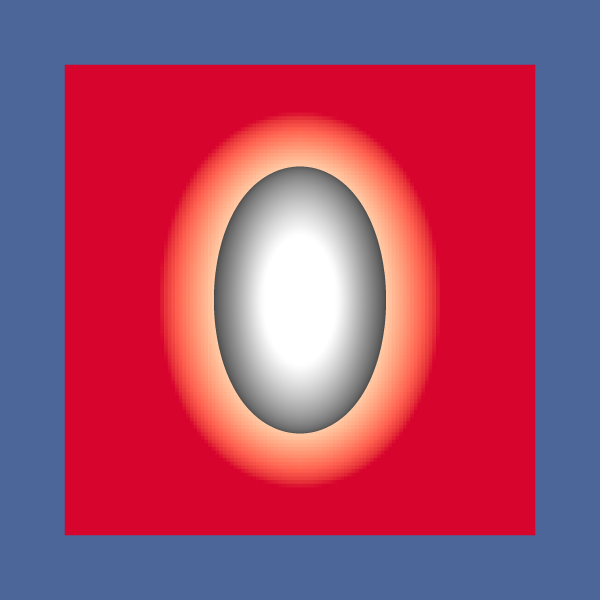}}
  \end{minipage}
  \caption{Final value of the 0.8-level-set}
  \label{fig:reinit_0.8LS_final}
\end{figure}

\begin{figure}[h!]
    \begin{minipage}[h!]{0.45\textwidth}
  \includegraphics[width=.9\textwidth]{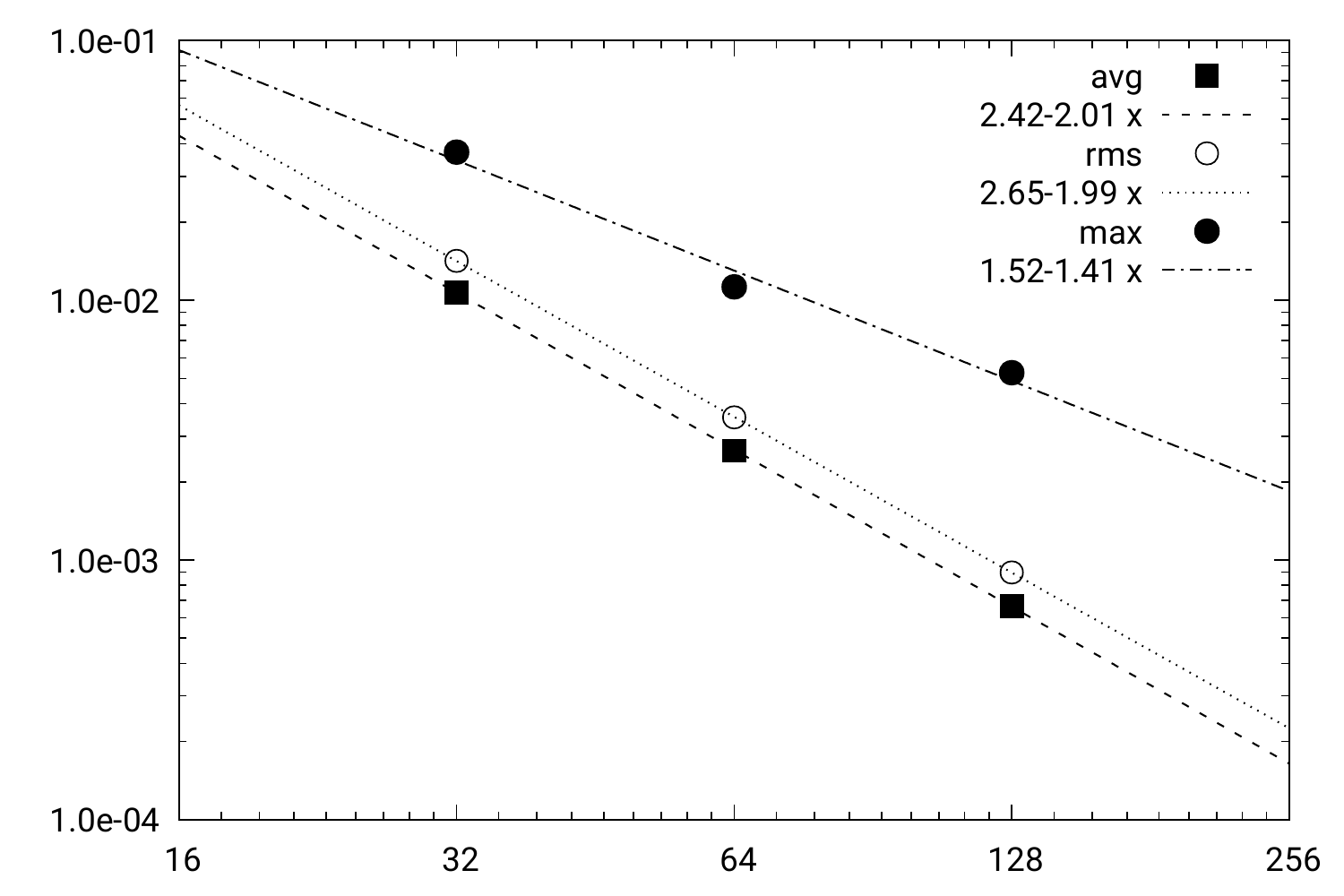}
  \subcaption{Error on the 0-level-set}
  \end{minipage}
  \begin{minipage}[h!]{0.45\textwidth}
  \includegraphics[width=.9\textwidth]{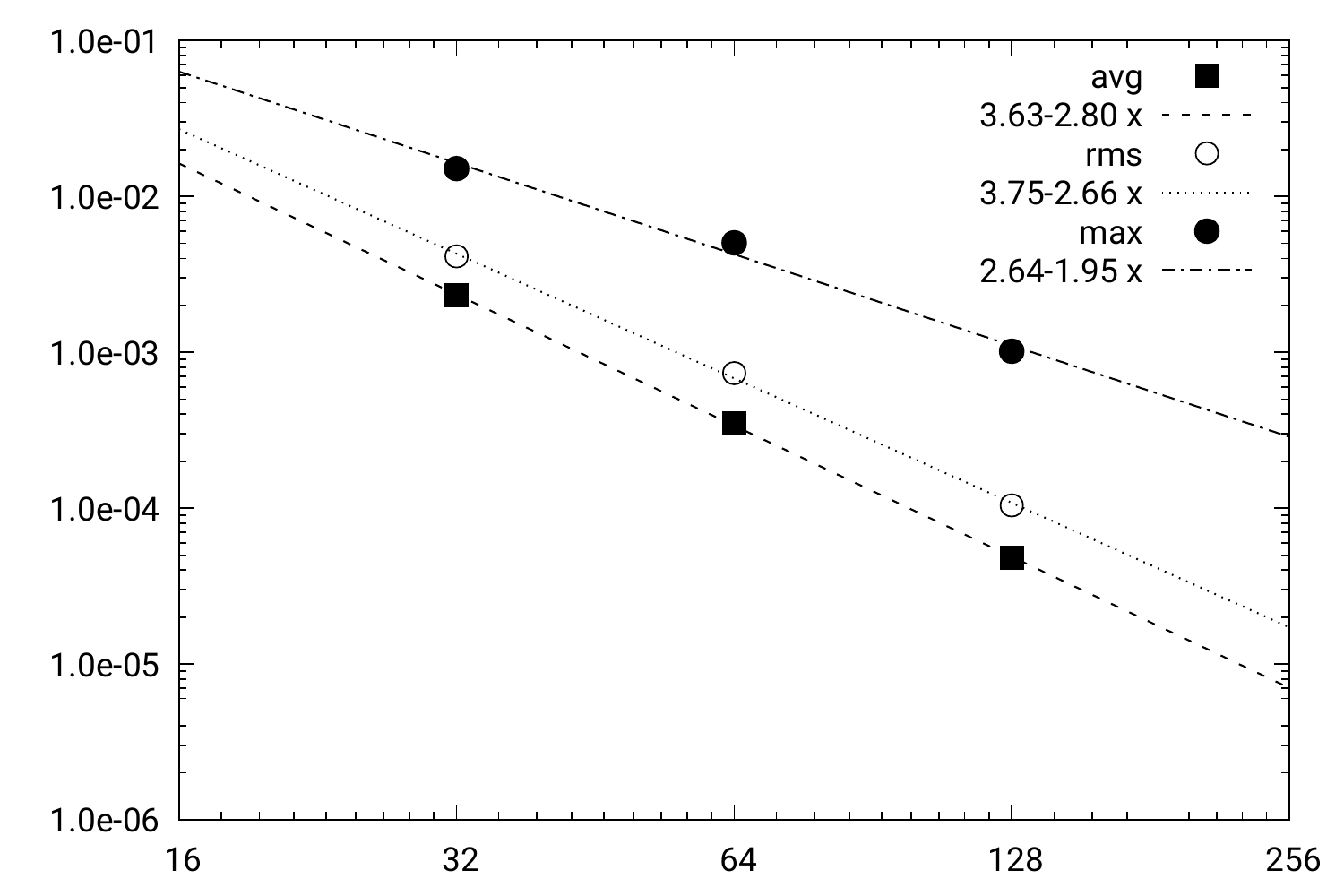}
  \subcaption{All cells}
  \end{minipage}
  \caption{Convergence results}
  \label{fig:reinit_convergence}
\end{figure}

The main idea behind this method derived from Russo \& Smereka
\cite{russo_remark_2000} is to have a subcell-accurate method in interfacial
cells and a simple spatial discretization operator elsewhere, for instance an
Essentially Non-Oscillatory (ENO) scheme
\begin{equation}
D_{x}^+\phi_{ijk} = \dfrac{\phi_{i+1jk} - \phi_{ijk}}{\Delta x} - 
\dfrac{\Delta}
{2}\minmod(D_{xx}\phi_{ijk},D_{xx}\phi_{i+1jk})
\label{eq:ENO_plus}
\end{equation}
\begin{equation}
D_{x}^-\phi_{ijk} = \dfrac{\phi_{ijk} - \phi_{i-1jk}}{\Delta x} - 
\dfrac{\Delta}
{2}\minmod(D_{xx}\phi_{ijk},D_{xx}\phi_{i-1jk})
\label{eq:ENO_minus}
\end{equation}
where
\begin{equation*}
D_{xx}\phi_{ijk} = \dfrac{\phi_{i-1jk} - 2\phi_{ijk} + \phi_{i+1jk}}
{\Delta x^2}
\end{equation*}
and
\begin{equation*}
\minmod(\alpha,\beta) = \begin{aligned}&\text{if } (\alpha\beta>0)\left\{
\begin{aligned}
&\text{if } |\alpha|<|\beta| \text{ , }\alpha\\
&\text{else } \beta 
\end{aligned}\right.\\
&\text{else } 0
\end{aligned}
\end{equation*}
we then define a Hamiltonian $H_G$ such that:
\begin{equation}
H_G(\boldsymbol{a},\boldsymbol{b}) = \left\{ \begin{aligned}
\sqrt{\sum_{i=1}^{d}\max((a_{i}^-)^2,(b_{i}^+)^2)} \text {\hspace{0.5cm} if } sgn
(\phi^0) \geq 0\\
\sqrt{\sum_{i=1}^{d}\max((a_{i}^+)^2,(b_{i}^-)^2)} \text {\hspace{0.5cm} if } sgn
(\phi^0) < 0
\end{aligned}
\right.
\end{equation}
where $d$ is the number of dimensions of the problem considered, $a_{i}^{+} =
\max(a_{i},0)$, $a_{i}^{-}=\min(a_{i},0)$ and $\boldsymbol{a},\boldsymbol{b}$ are
vectors such that
\begin{eqnarray}
\boldsymbol{a}&= (a_{i}) = (D^{+}_{i}\phi) \text{ , } i = \{x,y,z\}\\
\boldsymbol{b}&= (b_{i}) = (D^{-}_{i}\phi) \text{ , } i = \{x,y,z\}.
\end{eqnarray}
Thus, \cref{eq:LS_redist} becomes:
\begin{equation}
\phi_{\tau} + sign(\phi^0)[H_G(a,b)] = 0
\end{equation}
the ENO scheme \cref{eq:ENO_plus,eq:ENO_minus} is modified in cells where the
interface is located to limit the displacement of the
0-level-set. A quadratic ENO polynomial interpolation gives:
\begin{equation}
D_x^+\phi_{ijk} = \dfrac{0-\phi_{ijk}}{\Delta x^+} - \dfrac{\Delta x^+}{2}
\minmod
(D_{xx}\phi_{ijk},D_{xx}\phi_{i+1jk})
\end{equation}
and
\begin{equation}
\Delta x^+ = \left\{ \begin{aligned}
\Delta x &\text{ }\left( \dfrac{\phi^0_{i,j}-\phi^0_{i+1jk}-sgn(\phi^0_
{ijk}-\phi^0_{i+1jk})\sqrt{D}}{\phi^{0}_{xx}}\right) \text{ if } \left| \phi^0_
{xx}\right| >\epsilon \\
\Delta x &\text{ } \dfrac{\phi^0_{ijk}}{\phi^0_{ijk}-\phi^0_{i+1jk}} \text{
else.}\\
\end{aligned}
\right.
\end{equation}
with
\begin{eqnarray*}
\phi_{xx}^0 &=& \minmod(\phi^0_{i-1jk}-2\phi^0_{ijk}+\phi^0_{i+1jk}
\text{\hspace{0.2cm},\hspace{0.2cm}}
\phi^0_{ijk}-2\phi^0_{i+1jk}+\phi^0_{i+2jk}) \\
D &=& \left( \phi^0_{xx}/2  - \phi_{ijk}^0 - \phi_{i+1jk} \right)^2  - 4\phi_
{ijk}^0\phi_{i+1jk}^0
\end{eqnarray*}
$D^-_x\phi_{ijk}$ is modified in a similar fashion, see \cite{Min2010} for
details. One should note that for a smooth interface, without kinks, this method
is third-order accurate, whereas in the interfacial cells the order of accuracy is
reduced and falls between 1 and 2. We validated our method with a 3D case
adapted from \cite{russo_remark_2000} of a perturbed distance field to an
ellispoid:
\begin{equation}
\phi(x,y,z,t=0) = f(x,y,z) \times g(x,y,z)
\end{equation}
with $g$ real distance, and $f$ a perturbation function such that:
\begin{equation}
f(x,y) = \epsilon  + (x - x_0)^2 +(y - y_0)^2 + (z - z_0)^2
\end{equation}
with $x_0 = 3.5$, $y_0 = 2.$, $z_0 = 1.$, $\epsilon = 0.1$. Results are shown on
\cref{fig:reinit_0LS,fig:reinit_0.8LS,fig:reinit_0.8LS_final,fig:reinit_0LS_final} 
where we plot the isosurface a certain level-set with a slice view of the
distance function before and after the reinitialization. We show here that
we have extended Min's method to 3D calculations\footnote{The associated code for redistanciation is available at: \href{http://basilisk.fr/sandbox/alimare/LS_reinit.h}{LS\_reinit.h} and the associated ellipsoid redistanciation can be found at: \href{http://basilisk.fr/sandbox/alimare/1_test_cases/distanceToEllipsoid.c}{distanceToEllispoid.c}.} and obtain an order of accuracy
of 2 on \cref{fig:reinit_convergence}. 

For the discretization in time we use the TVD RK3 of \cite{Shu1988}. In Min
\cite{Min2010}, the author demonstrated that the fastest and most accurate
method is to use a Gauss-Seidel iteration with a fast-sweeping method
\cite{Tsai2003}. The raster-scan visiting algorithm associated (loops going
from $N_{x}$ to 1) would require specific cache construction to work with our
\textit{foreach()} iterators on adaptive grids which would probably compensate
the gains associated with this method.

\end{document}